\documentclass[12pt,a4paper]{amsart}
\usepackage[latin1]{inputenc}
\usepackage{latexsym}
\usepackage{color,graphicx,shortvrb}
\usepackage{amsmath, amssymb}
\usepackage{amsfonts}
\usepackage[colorlinks, bookmarks=true]{hyperref}

\usepackage[active]{srcltx}

\newtheorem{theorem}{Theorem}[section]
\newtheorem{lemma}[theorem]{Lemma}
\newtheorem{proposition}[theorem]{Proposition}
\newtheorem{corollary}[theorem]{Corollary}

\theoremstyle{definition}
\newtheorem{definition}[theorem]{Definition}
\newtheorem{remark}[theorem]{Remark}

\newtheorem{example}[theorem]{Example}

\newcommand{\field}[1] {\mathbb{#1}}
\newcommand{\R}{\field{R}}

\newcommand{\N}{\mathbb{N}}
\newcommand{\C}{\mathbb{C}}

\newcommand{\D}{\mathbb{D}}
\newcommand{\T}{\mathbb{T}}
\newcommand{\one}{\mathbf{1}}
\newcommand{\spn}{\mathbf{span}}
\newcommand{\rank}{\mathbf{rank}}

\newcommand{\se}{\mathcal{E}}
\newcommand{\is}{\mathcal{S}}
\newcommand{\ce}{\mathcal{S}_{\se}}
\newcommand{\moco}{{\mathrm{C}}}

\title{Shapiro's Theorem for subspaces}
\author{J. M. Almira and T. Oikhberg}

\begin{document}


\baselineskip=16pt

\numberwithin{equation}{section}

\maketitle \markboth{Shapiro's Theorem for subspaces}{J. M. Almira, T. Oikhberg}

\begin{abstract}
In their previous paper \cite{almira_oikhberg}, the authors investigated the existence of an element
$x$ of a quasi-Banach space $X$ whose errors of best approximation by a given
approximation scheme $(A_n)$ (defined by
$E(x,A_n) = \inf_{a \in A_n} \|x - a_n\|$) decay arbitrarily slowly. In this
work, we consider the question of whether $x$ witnessing the slowness rate
of approximation can be selected in a prescribed subspace of $X$.
In many particular cases, the answer turns out to be positive.
\end{abstract}

\section{Introduction, and an outline of the paper}
Let $(X,\|\cdot\|)$ be a quasi-Banach space, and let
$A_0\subset A_1\subset\ldots \subset A_n\subset\ldots \subset X$
be an infinite chain of subsets of $X$, where all inclusions are
strict. We say that $(X,\{A_n\})$ is an {\it approximation scheme}
(or that $(A_n)$ is an approximation scheme in $X$) if:
\begin{itemize}
\item[$(i)$] There exists a map $K:\mathbb{N}\to\mathbb{N}$ such that
$K(n)\geq n$ and $A_n+A_n\subseteq A_{K(n)}$ for all $n\in\mathbb{N}$.

\item[$(ii)$] $\lambda A_n\subset A_n$ for all $n\in\mathbb{N}$ and
all scalars $\lambda$.

\item[$(iii)$] $\bigcup_{n\in\mathbb{N}}A_n$ is a dense subset of $X$
\end{itemize}
An approximation scheme is called {\it non-trivial} if $X \neq \cup_n \overline{A_n}$.

Many problems in approximation theory can be described using approximation schemes.
We say that an approximation scheme is {\it linear} if the
sets $A_n$ are linear subspaces of $X$. In this setting, we can take $K(n) = n$.
Linear approximation schemes arise, for instance, in problems of approximation
of functions by polynomials of prescribed degree. Non-linear schemes arise,
for instance, in the context of the so-called adaptive approximation by
elements of a dictionary (see Definition~\ref{def_dictionary} below).
R.~DeVore's survey paper
\cite{devorenonlinear} provides a good introduction into
adaptive approximation and its advantages.

Approximation schemes were introduced by Butzer and Scherer in 1968
\cite{butzer_scherer} and, independently,  by Y. Brudnyi and N. Kruglyak under the
name of ``approximation families'' in 1978 \cite{brukru}, and popularized by Pietsch
in his seminal paper of 1981 \cite{Pie}, where the approximation spaces
${\mathbf{A}}_p^r(X,{A_n})=\{x\in X: \|x\|_{A_p^r}=
\|\{E(x,A_n)\}_{n=0}^\infty\|_{\ell_{p,r}}<\infty\}$ were studied.
Here, $$\ell_{p,r}=\{\{a_n\}\in\ell_\infty: \|\{a_n\}\|_{p,r}
=\left[\sum_{n=1}^\infty n^{rp-1}(a_n^*)^p\right]^\frac{1}{p}<\infty\}$$
denotes the so called Lorentz sequence space, and $E(x,A_n)=\inf_{a\in A_n}\|x-a\|_X$.
In \cite{Pie}, it was also proved that
${\mathbf{A}}_p^r(X,{A_n})\hookrightarrow {\mathbf{A}}_q^s(X,{A_n})$
holds whenever $r>s>0$, or $r=s$ and $p<q$ (in other words, the approximation
spaces form a scale).

In the context of approximation of functions by polynomials, the classical
theorems of Bernstein and Jackson (see e.g.~\cite[Section 7]{devore}) indicate
a strong connection between the membership of a function $f$ in a space
${\mathbf{A}}_p^r(X,{A_n})$, and the degree of smoothness of $f$.
For this reason, the spaces ${\mathbf{A}}_p^r(X,{A_n})$ are often
referred to as ``generalized smoothness spaces'' (see, for example,
\cite{devorenonlinear}, \cite{devore}, \cite{PieID}). Thus, we can
view the rate of decrease of a sequence $(E(x,A_n))$ as reflecting the
``smoothness'' of $x$.

To proceed further, we fix some notation. We write $\{\varepsilon_i\} \searrow 0$
to indicate that the sequence $\varepsilon_1 \geq \varepsilon_2 \geq \ldots \geq 0$
satisfies $\lim_i \varepsilon_i = 0$. For a quasi-normed space $X$, we denote by $B(X)$
and $S(X)$ its closed unit ball and unit sphere, respectively. That is,
$S(X) = \{x \in X : \|x\| = 1\}$, and $B(X) = \{x \in X : \|x\| \leq 1\}$.
We use the notation $\mathbf{B}(X,Y)$ for the space of bounded linear operators $T:X\to Y$, with the usual convention $\mathbf{B}(X)=\mathbf{B}(X,X)$.
If $X$ is a quasi-Banach space, $x \in X$, and $A \subset X$, we define
the {\it best approximation error} by
$E(x,A)_X = dist(x,A)_X = \inf_{a \in A} \|x-a\|$.
When there is no confusion as to the ambient space $X$ and its (quasi-)norm,
we simply use the notation $E(x,A)$. If $B$ and $A$ are two subsets of $X$,
we set $E(B,A) = \sup_{b \in B} E(b,A)$ (note that $E(B,A)$ may be different
from $E(A,B)$).

The results described below have their origin in
the classical Lethargy Theorem by
S.N. Bernstein \cite{bernsteininverso}, stating that, for any linear
approximation scheme $(A_n)$ in a Banach space $X$ ,
if $\dim A_n<\infty$ for all $n$ and $\{\varepsilon_n\}$ is a non-increasing
sequence of positive numbers, $\{\varepsilon_n\}\in c_0$, there exists
$x \in X$ such that $E(x,A_n)=\varepsilon_n$ for all $n\in \N$.
Bernstein's proof was based on a compactness argument, and only
works if $\dim A_n<\infty$ for all $n$. In 1964 H.S. Shapiro
\cite{shapiro} used Baire Category Theorem and Riesz's Lemma (on the existence of
almost orthogonal elements to any closed linear subspace $Y$ of a Banach space $X$)
to prove that, for any sequence $A_1 \subsetneq A_2 \subsetneq \ldots \subsetneq X$ of
closed (not necessarily finite dimensional) subspaces of a Banach space $X$,
and any sequence $\{\varepsilon_{n}\}\searrow 0$, there exists an $x\in X$ such that
$E(x,A_{n})\neq\mathbf{O}(\varepsilon_{n})$.  This result was strengthened by
Tjuriemskih \cite{tjuriemskih1}, who, under the very same conditions of Shapiro's
Theorem, proved the existence of  $x\in X$ such that $E(x,A_{n})\geq \varepsilon_{n}$,
$n=0,1,2,\ldots$. Later, Borodin \cite{borodin} gave an elementary proof of this
result. He also proved that, for arbitrary infinite dimensional Banach spaces $X$,
and for any sequence $\{\varepsilon_n\}\searrow 0$ satisfying
$\varepsilon_n>\sum_{k=n+1}^\infty\varepsilon_k$, $n=0,1,2,\ldots$, there exists
$x\in X$ such that $E(x,X_{n})= \varepsilon_{n}$, $n=0,1,2,\ldots$.

Motivated by these results, in \cite{almira_oikhberg} the authors
gave several characterizations of the approximation schemes $(X,\{A_n\})$ with
the property that for every non-increasing sequence
$\{\varepsilon_n\}\searrow 0$ there exists an element $x\in X$ such that
$E(x,A_n) \neq \mathbf{O}(\varepsilon_n)$. In this case we say
that $(X,\{A_n\})$ (or simply $(A_n)$) {\it satisfies Shapiro's Theorem}.
We established the following characterization of
approximation schemes satisfying Shapiro's Theorem
(see \cite[Theorem 2.2, Corollary 3.6]{almira_oikhberg}):

\begin{theorem}\label{teo_introd}
Let $X$ be a quasi-Banach space.
For any approximation scheme $(X,\{A_n\})$, the following are equivalent:
\begin{itemize}
\item[$(a)$]  The approximation scheme $(X,\{A_n\})$ satisfies Shapiro's Theorem.

\item[$(b)$] There exists a constant $c>0$ and an infinite set $\mathbb{N}_{0}\subseteq\mathbb{N}$ such that for all $n\in\mathbb{N}_{0}$, there
exists some $x_{n}\in X\setminus \overline{A_{n}}$ which satisfies $E(x_{n},A_{n})\leq cE(x_{n},A_{K(n)}).$

\item[$(c)$]
There is no decreasing sequence $\{\varepsilon_n\}\searrow 0$ such that
$E(x,A_n)\leq \varepsilon_n\|x\|$  for all $x\in X$ and $n\in \N$.

\item[$(d)$] $E(S(X),A_n)=1$, $n=0,1,2,\ldots$.

\item[$(e)$] There exists $c>0$ such that $E(S(X),A_n)\geq c$, $n=0,1,2,\ldots$.
\end{itemize}
Moreover, if $X$ is a Banach space, then all these conditions are equivalent to:
\begin{itemize}
\item[$(f)$] For every non-decreasing sequence $\{\varepsilon_n\}_{n=0}^{\infty}\searrow 0$ there exists an element $x\in X$ such that $E(x,A_n)\geq \varepsilon_n$ for all $n\in\mathbb{N}$.
\end{itemize}
\end{theorem}

Riesz's Lemma claims that condition $E(S(X),A_n)=1$, appearing at item
$(d)$ above, holds whenever $X$ is a Banach space and $A_n$ is a closed linear
subspace of $X$. Therefore, any  non-trivial linear approximation scheme $(A_n)$ in
a Banach space $X$ satisfies Shapiro's Theorem. Thus, Theorem \ref{teo_introd}
generalizes Shapiro's original result \cite{shapiro}.


In this paper, we consider Shapiro's theorem in the setting of
{\it constrained approximation}. To the best of our knowledge,
``constrained'' versions of lethargy theorems have never been studied.
Indeed, a search of Mathscinet for the years from 2000
to 2010 yielded $122$ items with primary AMS classification $41A29$
(approximation with constraints), none of them dealing with lethargy
problems. To fill this gap, in this paper we investigate the following
``restricted'' version of Shapiro's Theorem.

\begin{definition}\label{subspaces}
Suppose $Y$ is a linear subspace of a quasi-Banach space $X$. We say that $Y$
{\it satisfies Shapiro's Theorem} with respect to the approximation scheme
$(X,\{A_n\})$ if, for any $\{\varepsilon_n\} \searrow 0$,
there exists $y \in Y$ such that $E(y, A_n)_X \neq {\mathbf{O}}(\varepsilon_n)$.
\end{definition}


By default, we view $Y$ as a space, equipped with its own quasi-norm, and
embedded continuously into $X$. If, in addition, $Y$ is a closed subspace
of $X$, Open Mapping Theorem (see \cite[Corollary 1.5]{kalton}) shows that
the norms $\|\cdot\|_X$ and $\|\cdot\|_Y$ are equivalent on $Y$.

This paper is organized as follows.
We start by giving a general description
of subspaces satisfying Shapiro's Theorem (Section~\ref{criteria}).
One of our main tools is the notion of $Y$ being ``far'' from
an approximation scheme $(A_n)$ (Definition~\ref{def_far}).
We show that if $Y$ satisfies Shapiro's Theorem
relative to the approximation scheme $(A_n)$, then $Y$ is $c$-far from $(A_n)$
for a certain positive constant $c>0$. If $Y$ is a closed subspace of $X$, the
converse is also true (Theorem \ref{far_shap}).
We use this characterization to prove
that, if $(X,\{A_n\})$ satisfies Shapiro's Theorem,
then all finite codimensional subspaces
of $X$ satisfy Shapiro's Theorem relative to $(A_n)$ (Theorem \ref{fin_codim}).
On the other hand, ``small'' subspaces (for instance, subspaces of $X$ of
countable algebraic dimension) fail Shapiro's
Theorem (Corollary \ref{countable}).
We end Section~\ref{criteria} by noting a link between the notion of being far,
and a generalized version of the classical theorem of Jackson, connecting
the rate of approximation of a function with its degree of smoothness
(see Proposition~\ref{far_jackson}, and the remarks preceding it).

Section~\ref{complemsub} deals with the case when there exists a bounded
projection $P$ from $X$ onto $Y$. Theorem \ref{subsp_shapiro} gives
several criteria for $Y$ to satisfy Shapiro's Theorem relative to $(P(A_n))$.
It also shows that, if $Y$ satisfies Shapiro's Theorem relative to $(P(A_n))$,
then $Y$ satisfies Shapiro's Theorem relative to $(A_n)$.
Theorem \ref{1_far_for_stable_schemes} shows that, if
$Y$ has finite codimension, and $(X,\{A_n\})$ satisfies Shapiro's Theorem, then
$Y$ satisfies Shapiro's Theorem relative to $(P(A_n))$.
Along the way, we prove that an interesting stability result:
if an approximation schemes $(A_n)$ in $X$ satisfies Shapiro's Theorem, and
$F$ is a finite dimensional subspace of $X$, then the scheme $(A_n+F)$
satisfies Shapiro's Theorem, too (Theorem \ref{stability}).

Section~\ref{compact_schemes} is devoted to boundedly compact approximation schemes
$(X,\{A_n\})$ (that is, $B(X) \cap A_n$ is relatively compact in $X$, for every $n$).
In this case any infinite dimensional closed subspace of $X$ satisfies Shapiro's Theorem
(Theorem \ref{compact_shap}). If, furthermore, the sets $A_n$ are linear finite
dimensional subspaces of $X$, then, for any infinite
dimensional closed subspace $Y$ of $X$, and any sequence
$\{\varepsilon_n\}\in c_0$, there is an element $y\in Y$ such that
$E(y,A_n)\geq \varepsilon_n$ for all $n$ (Theorem \ref{bern_subspace}).

In Section \ref{sec_compact} we study the subspaces $Y$ compactly embedded into $X$.
In this case, $Y$ cannot be far from any approximation scheme
$(A_n)$ (Theorem~\ref{compact_subspaces}), hence it fails Shapiro's Theorem.
If $(A_n)$ is boundedly compact, then the spaces $Y$ failing Shapiro's
Theorem are precisely those that are included into a compactly embedded subspace
$Z$ of $X$ (Theorem~\ref{compact_embed}). Several examples of compactly embedded
subspaces are provided.

Section~\ref{1-far} describes approximation schemes
$(X,\{A_n\})$ with the property that all finite
codimensional subspaces $Y$ of $X$ are $1$-far from  $(A_n)$.
The main characterization is given by Theorem~\ref{basis_far}.
As an aid of our investigation, we introduce and study the
Defining Subspace Property of Banach spaces.

Finally, in Section~\ref{examples}, we exhibit several additional examples
subspaces (arising from harmonic analysis) which satisfy Shapiro's Theorem.

Note that we encounter several instances of continuous functions on $[a,b]$,
analytic on $(a,b)$, which are ``poorly approximable''
(Corollaries \ref{suavidad}, \ref{cor_gen_haar}). This illustrates the thesis
that the smoothness conditions guaranteeing that a function is ``well approximable''
must be ``global.'' The failure of smoothness at endpoints may result in an
arbitrarily slow rate of approximation.


Throughout the paper, we freely use standard functional analysis facts and notation.
Recall that, if $\|\cdot\|$ is a quasi-norm on the vector space $X$, then there is a constant $C_X\geq 1$
such that the inequality $\|x_1 + x_2\| \leq C_X (\|x_1\| + \|x_2\|)$
holds for any $x_1, x_2 \in X$ (the usual triangle inequality occurs when $C_X = 1$).
The space $X$ is called {\it $p$-convex} ($0 < p \leq 1$) if
$\|x_1 + x_2\|^p \leq \|x_1\|^p + \|x_2\|^p$ for any $x_1, x_2 \in X$
(any normed space is $1$-convex). The classical Aoki-Rolewicz theorem states that
every quasi-Banach space has an equivalent $p$-convex norm, for some $p$ \cite{kalton}.
If $A$ is a subset of the quasi-normed space $X$, we denote by $\spn[A]$ the algebraic
linear span of $A$, and by 
$\overline{A}$ its quasi-norm closure.

\section{Criteria for Shapiro's Theorem}\label{criteria}

In this section we investigate general properties of subspaces
satisfying Shapiro's Theorem.
One of our main tools is the notion of a subspace being ``far''
from an approximation scheme.

\begin{definition}\label{def_far}
Suppose $(A_n)$ is an approximation scheme in $X$, and a quasi-normed space $Y$ is
embedded continuously into $X$.
We say that $Y$ is {\it $c$-far} from $(A_n)$ if $E(S(Y), A_n)$ $\geq c$
for every $n$. 
The subspace $Y$ is said to be {\it far} from $(A_n)$ if it is $c$-far
from $(A_n)$ for some $c>0$. We say that $Y$ is {\it not far} from $(A_n)$
if there is no $c > 0$ with the property that $Y$ is $c$-far from $(A_n)$.
\end{definition}

Theorem~\ref{far_shap} shows that, if $Y$ satisfies Shapiro's Theorem relative
to $(A_n)$, then it is far from $(A_n)$. The converse is true if $Y$ is closed.
We then prove that Shapiro's Theorem and ``farness'' are stable under isomorphisms
(see e.g. Proposition \ref{shapiro_pert}), but not under contractive embeddings
(Proposition \ref{subspace_shapiro}). We prove
that, in some cases, ``large'' (for instance, finite codimensional) subspaces
of $X$ must be far from approximation schemes (Proposition \ref{fin_codim_p}),
and must satisfy Shapiro's Theorem (Theorem \ref{fin_codim}). On the contrary,
``small'' subspaces fail Shapiro's Theorem (Corollary \ref{countable}).
Finally, we note that ``farness'' can be viewed as a generalization of
classical results of Jackson on the approximation of smooth functions
(Proposition \ref{far_jackson}).

\begin{theorem}\label{far_shap}
Suppose $(A_n)$ is an approximation scheme in $X$, and a quasi-normed space $Y$ is
embedded continuously into $X$.
\begin{enumerate}
\item
If $Y$ satisfies Shapiro's Theorem relative to $(A_n)$, then it is far from $(A_n)$.
\item
Conversely, every closed subspace of $X$ which is far from $(A_n)$, satisfies
Shapiro's Theorem relative to $(A_n)$.
\end{enumerate}
\end{theorem}

This theorem states that a subspace $Y$, satisfying Shapiro's Theorem
relative to $(A_n)$, must be $c$-far from $(A_n)$, for some $c \in (0,1]$.
By Remark \ref{not-1-far}, this $c$ can be arbitrarily close to $1$.
In Section \ref{1-far}, we investigate the ``extreme case'' of subspaces which
are $1$-far from approximation schemes.

\begin{proof}
(1)
Suppose first that $Y$ is not far from $(A_i)$, and show the failure of Shapiro's
Theorem. Indeed, in this case, there exists a sequence $0 = i_1 < i_2 < \ldots$,
such that $E(S(Y), A_{i_k}) \leq 1/k$ for $k \in \N$. Define $\varepsilon_i = 1/k$
for $i_k \leq i < i_{k+1}$. Then $E(y,A_i) \leq \varepsilon_i \|y\|$ for any
$y \in Y$. In other words, the sequence $\{\varepsilon_i\} \searrow 0$
witnesses the failure of Shapiro's Theorem.

(2)
Now suppose $Y$ is a closed subspace of $X$ (equipped with the norm
inherited from $X$), which is far from $(A_i)$.
Renorming if necessary, we can assume that $X$ is $p$-convex.
Find $c \in (0,1)$ such that, for every $i$, there exists $y \in Y$ satisfying
$c < E(y,A_i) \leq \|y\| < 1$. For a given sequence $\{\varepsilon_n\} \searrow 0$
of positive numbers, let us see that we can find a sequence $0 = i_0 < i_1 < \ldots$, and $y \in Y$,
such that $E(y, A_{i_j}) \geq 2^{j-1} \varepsilon_{i_j}$ for every $j \geq 1$.

Define the sequence $(i_j)$ recursively. Set $i_0 = 0$. Pick $i_1 \in \N$
such that $\varepsilon_{i_1} < c \varepsilon_0/8^{1/p}$.
Find $y_1 \in A_{i_1}$ with $c < \|y_1\| < 1$.

Suppose
$i_0 < \ldots < i_{j-1}$ have already been selected. Let $s_j = K^{j}(i_{j-1})$,
where $K^j = K \circ \ldots \circ K$ ($j$ times).
Pick $i_j > s_j$ in such a way that (i) there exists $y_j \in A_{i_j}$
satisfying $\|y_j\| < 1$ and $E(y_j, A_{K(s_j)}) > c$, and
(ii) $\varepsilon_{i_j} < c \varepsilon_{i_{j-1}}/8^{1/p}$.

For $j \geq 1$ let $\alpha_j = 2^{j/p} c^{-1} \varepsilon_{i_{j-1}}$.
Then, for $m > j$, $\alpha_m < c \alpha_j/4^{(m-j)/p}$.
Set $y = \sum_{j=1}^\infty \alpha_j y_j$
(the series converges, since $\sum_j \alpha_j^p < \infty$).
Then, for any $j$,
\begin{align*}
E(y, A_{i_{j-1}})^p
&
\geq
E(\sum_{k=0}^j \alpha_k y_k , A_{i_{j-1}})^p -
\sum_{k=j+1}^\infty \alpha_k^p \geq
E(\alpha_j y_j , A_{s_j})^p - \sum_{k=j+1}^\infty \alpha_k^p
\\
&
\geq
\alpha_j^p c^p - \sum_{k=j+1}^\infty \alpha_j^p c^p 4^{j-k} \geq
\frac{\alpha_j^p c^p}{2} > 2^{(j-1)p} \varepsilon_{i_{j-1}}^p .
\end{align*}
\end{proof}
\begin{remark}\label{part_3_theorem21}
The hypothesis of $Y$ being closed in $X$ can not be omitted from
Theorem \ref{far_shap}$(2)$.
More precisely, there exists a continuous embedding of a Banach space $Y$ to
a Banach space $X$, and an approximation scheme $({\mathcal{A}}_n)$ in $X$, such that
$Y$ fails Shapiro's Theorem with respect to $({\mathcal{A}}_n)$, but
$E(S(X)\cap Y, A_n) = 1$ for every $n$.  For instance,
$X=C[0,2\pi]$. For $n \in \N$, let $\mathcal{A}_n$ denote
the space of algebraic polynomials of degree less than $n$. For $1 \leq r < \infty$,
$$
Y={\mathbf{A}}_r^r(C[0,2\pi],\{\mathcal{A}_n\}_{n=0}^\infty)=
\Big\{f\in C[0,2\pi]:
\|f\|_Y:= \big(\sum_n E(f,\mathcal{A}_n)^r\big)^{1/r} < \infty\Big\}
$$
is an infinite dimensional Banach space \cite[Section 3]{almiraluther2}.
Furthermore, $Y$ is continuously embedded into $X$.
As the sequence $(E(f,\mathcal{A}_n))_n$ is non-increasing,
$E(f,\mathcal{A}_n) \leq (n+1)^{-r} \|f\|_Y$ for any $f \in Y$.
Thus, $Y$ fails Shapiro's Theorem for $(\mathcal{A}_n)$. Moreover,
$(C[0,2\pi],\{\mathcal{A}_n\}_{n=0}^\infty)$ is a non-trivial
linear approximation scheme, hence it satisfies Shapiro's Theorem.
To show that $E(S(X)\cap Y, \mathcal{A}_n) = 1$, let $h(t) = \cos nt$. Then
$h \in S(X) \cap Y$. By Chebyshev Alternation Theorem
(see e.g. \cite[Section 3.5]{devore}), $E(h, \mathcal{A}_n) = 1$.
\end{remark}

\begin{remark}\label{not-1-far}
For any $c \in (0,1)$ one can find a linear approximation scheme $(A_n)$ in
$\ell_2$, and a closed subspace $Y$, which is $c$-far from $(A_n)$, but not $c_1$-far
if $c_1 > c$. Indeed, denote the canonical basis for $\ell_2$ by $(e_i)$.
For $i \in \N$ let $f_i = \sqrt{1-c^2} \, e_{2i} + c e_{2i-1}$. For $n \in \N$,
let $A_n$ be the closed linear span of the vectors $e_j$, where $j$ is either
even, or does not exceed $2n-2$. Let $Y$ be the closed linear span of
the vectors $f_i$. Clearly, for any $y \in S(Y)$ and $n \in \N$, $E(y,A_n) \leq c$.
Furthermore, $E(f_n,A_n) = c$ for every $n$.
\end{remark}



Next we show that subspaces satisfying Shapiro's Theorem are
stable under small perturbations.

\begin{lemma}\label{far_pert}
Suppose $Y$ and $Z$ are subspaces of a $p$-convex quasi-Banach space $X$,
equipped with the norm inherited from $X$.
Suppose, furthermore, that $Y$ is $c$-far from an approximation scheme $(A_i)$,
and that $E(S(Y),Z)<c$.  Then $Z$ is $c_1$-far from $(A_i)$, with
$c_1=\left(\frac{c^p-E(S(Y),Z)^p}{1+c^p}\right)^{1/p}$.
\end{lemma}

An application of Aoki-Rolewicz theorem then yields the following.

\begin{proposition}\label{shapiro_pert}
Suppose $Y$ is a closed subspace of a quasi-Banach space $X$, satisfying
Shapiro's Theorem relative to an approximation scheme $(A_i)$. Then
there exists $\delta > 0$ such that any closed subspace $Z$ of $X$,
with the property that $E(S(Y),Z) < \delta$, also satisfies Shapiro's Theorem
relative to $(A_i)$.
\end{proposition}
\begin{proof}[Proof of Lemma~\ref{far_pert}]
For the sake of brevity, set $E=E(S(Y),Z)$. Then for any $\lambda>0$ with
$\lambda^p \in (0,c^p-E^p)$ there exist $\alpha,\beta\in (E,c)$ such that $\beta>\alpha$ and $\beta^p-\alpha^p>\lambda^p$. Now, $\beta<c$ implies that, for each $i\in\N$ there exists $y \in S(Y)$ such that $E(y,A_i) > \beta$. On the other hand, $\alpha>E$ implies that there exists $w \in Z$ such that $\|y - w\| < \alpha$. Hence
\[
\beta^p<E(y,A_i)^p=\inf_{a_i\in A_i}\|y-w+w-a_i\|^p\leq \|y-w\|^p+E(w,A_i)^p\leq \alpha^p+E(w,A_i)^p.
\]
Moreover, $\|w\|^p\leq \|y\|^p+\|y-w\|^p\leq 1+\alpha^p\leq 1+c^p$. It follows that $E(\frac{w}{\|w\|},A_i)^p \geq \frac{\beta^p-\alpha^p}{1+\alpha^p}\geq \frac{\lambda^p}{1+c^p}$, and this holds for any $\lambda^p \in (0,c^p-E^p)$. Hence $E(S(Z),A_i)\geq \left(\frac{c^p-E(S(Y),Z)^p}{1+c^p}\right)^{1/p}$. This ends the proof.
\end{proof}


Intuitively, ``large'' subspaces of $X$ must be far from approximation schemes,
and must satisfy Shapiro's Theorem. Proposition \ref{fin_codim_p} and
Theorem \ref{fin_codim} prove these statements, in some cases.

\begin{proposition}\label{fin_codim_p}
Suppose $X$ is a $p$-convex quasi-Banach space ($p \in (0,1]$).
Consider an approximation scheme $(A_n)$ in $X$, satisfying Shapiro's Theorem,
and let $Y$ be a finite codimensional closed subspace of $X$. Then
$Y$ is $2^{-1/p}$-far from $(A_n)$.
\end{proposition}

Note that, if $Y$ is ``nicely complemented'' in $X$, the estimates of
Lemma \ref{far_pert} and Proposition \ref{fin_codim_p} can be improved
(Theorem \ref{1_far_for_stable_schemes}). Furthermore,
$1$-far subspaces are studied in Section \ref{1-far}.

For the proof of Proposition \ref{fin_codim_p} we need:

\begin{lemma}\label{dense}
Suppose $X$, $(A_n)$, $Y$, and $p$ are as in the statement of
Theorem~\ref{fin_codim_p}. Then for every $\delta > 0$ there exist
$a_1, \ldots, a_L \in \cup_n A_n$, such that for every $x \in B(X)$
there exist $\ell \in \{1, \ldots, L\}$ and $y \in 2^{1/p} (1+\delta) B(Y)$
satisfying $\|x - (a_\ell + y)\| < \delta$.
\end{lemma}

\begin{proof}
Find $c \in (0, ((1+\delta)^p - 1)^{1/p})$.
As $(\cup_n A_n) \cap B(X)$ is dense in $B(X)$, $q((\cup_n A_n) \cap B(X))$
is dense in $B(X/Y)$ (here, $q : X \to X/Y$ denotes the quotient map).
Thus, we can use that $\dim X/Y<\infty$ to find
$n \in \N$ and a $c/2$-net $(e_\ell)_{\ell=1}^L$ in
$B(X/Y)$, such that for any $\ell$ there exists $a_\ell \in A_n \cap B(X)$
with $q (a_\ell) = e_\ell$.

For any $x \in B(X)$ there exists $\ell \in \{1, \ldots, L\}$
such that  $\|q(x) - q(a_\ell)\| =E(x-a_\ell,Y) < c$. Hence there exists $y\in Y$ such that $\|x-a_\ell-y\|<c$. By the $p$-convexity of $X$,
\begin{align*}
\|y\|^p
&
=
\|y-(a_\ell-x)\|^p+\|a_\ell-x\|^p
\\
&
\leq \|y-(a_\ell-x)\|^p+\|a_\ell\|^p+ \|x\|^p\leq c^p+2\leq  (1+\delta)^p+1\leq 2(1+\delta)^p .
\end{align*}
This ends the proof.
\end{proof}

\begin{proof}[Proof of Proposition~\ref{fin_codim_p}]
Suppose our assumption is false. Then there exists $n \in \N$
such that $E(S(Y), A_n) < \gamma < 2^{-1/p}$.
Find $c \in (0, 2^{-1/p} - \gamma)$, in such a way that
$2 \gamma^p (1+c)^p + c^p < 1$.
By Lemma~\ref{dense}, there exists $m \in \N$ such that
for every $x \in B(X)$ there exists $a \in A_m$ and
$y \in 2^{1/p} (1+c) B(Y)$, satisfying $\|x - (a+y)\| < c$.
For $y$ as above, there exists $b \in A_n$ such that
$\|b - y\| < 2^{1/p} \gamma (1+c)$. Then $a+b \in A_N$,
where $N = K(\max\{n,m\})$. Furthermore, $x - (a+b) = (x - (a+y)) + (y - b)$,
hence, by our choice of $c$,
$$
\|x - (a+b)\|^p \leq \|x - (a+y)\|^p + \|y-b\|^p < 2 \gamma^p (1+c)^p + c^p .
$$
It follows that $E(S(X), A_N) < 1$, since $x\in B(X)$ was arbitrary. This contradicts
$(a) \Leftrightarrow (d)$ of Theorem~\ref{teo_introd}.
\end{proof}

Proposition~\ref{fin_codim_p} implies that, if $Y$ is a closed finite codimensional subspace
$X$, and the approximation scheme $(A_n)$ in $X$ satisfies Shapiro's Theorem, then $Y$
satisfies Shapiro's Theorem relative to $(A_n)$. In fact, a stronger result is true.

\begin{theorem}\label{fin_codim}
Let $(X,\{A_n\})$ be an approximation scheme satisfying Shapiro's Theorem and let  $Y$ be a finite codimensional subspace of $X$. Then:
\begin{itemize}
\item[$(1)$] $Y$ satisfies Shapiro's Theorem with respect to $(A_n)$.
\item[$(2)$] If $X$ is a Banach
space and $Y$ is closed,
then, for all non-increasing sequence of positive numbers $\{\varepsilon_n\}\in c_0$ there exists $y\in Y$ such that $E(y,A_n)\geq \varepsilon_n$ for all
$n\in\mathbb{N}$.
\end{itemize}
\end{theorem}

The result below (used to prove Theorem~\ref{fin_codim}) is of independent interest.

\begin{lemma}\label{sum_of_subsp}
Suppose $(A_i)$ is an approximation scheme in a quasi-Banach space $X$, and
$(Y_j)_{j \in I}$ is a finite or countable collection of subspaces of $X$,
each $Y_j$ failing Shapiro's Theorem relative to $(A_i)$.
Then $\spn[Y_j : j \in I]$ fails Shapiro's Theorem relative to $(A_i)$.
\end{lemma}

\begin{proof}
We present the proof for $I = \N$ (the finite
case is handled in a similar manner). 
As $X$ is a quasi-Banach space, there exists a constant $C_q \geq 1$
such that $\|x_1 + x_2\| \leq C_q(\|x_1\| + \|x_2\|)$ for any $x_1, x_2 \in X$.
A simple induction argument shows that, for any $x_1, \ldots, x_m \in X$,
we have $\|x_1 + \ldots + x_m\| \leq C_q^{m-1}(\|x_1\| + \ldots + \|x_m\|)$.
We shall write $K^j$ for $K \circ \ldots \circ K$ ($j$ times).

For every $k \in \N$ there exist a function $C_k : Y_k \to [0,\infty)$
and a sequence $\{\varepsilon_j^{(k)}\} \searrow 0$ such that
the inequality $E(y,A_j) \leq C_k(y) \varepsilon_j^{(k)}$ holds for every
$j \in \N$ and every $y \in Y_k$. Let $Y_m^\prime  = \spn[Y_i : i \leq m]$.
For any $y \in Y = \spn[Y_i : i \in \N]$, let $m(y)$ be the smallest $m \in \N$
for which $y \in Y_m^\prime$. Pick a representation $y = \sum_{\ell=1}^m y_\ell$
(with $m = m(y)$, and $y_\ell \in Y_\ell$), and set
$C(y) = \sum_{\ell=1}^m C_\ell(y_\ell)$. We shall construct a sequence
$\{\varepsilon_s\} \searrow 0$ such that $E(y,A_s) \leq C(y) \varepsilon_s$
for $s$ large enough. To this end, pick sequences
$$
0 = s_0 < n_1 \leq s_1 := K(n_1) < n_2 \leq s_2 := K^2(n_2) < n_3 \leq \ldots
$$
in such a way that $\varepsilon_{n_j}^{(k)} \leq (2 C_q)^{-j}$ for $1 \leq k \leq j$.
Then for any $y \in Y_m^\prime$ and $k \geq m$,
\begin{align*}
E(y, A_{s_k})
&
\leq
E(\sum_{\ell=1}^m y_\ell , A_{K^m(n_k)}) \leq
C_q^{m-1} \sum_{\ell=1}^m E(y_\ell, A_{n_k})
\\
&
\leq
C_q^{m-1} \sum_{\ell=1}^m C_\ell(y_\ell) \frac{1}{(2 C_q)^k} \leq
2^{-k} C(y) .
\end{align*}
Let $\varepsilon_s = 2^{-k}$ for $s_k \leq s < s_{k+1}$.
Then, for $y \in Y$ and $s \geq s_{m(y)}$,
$E(y,A_s) \leq E(y, A_{s_k}) \leq 2^{-k} = \varepsilon_s$.
\end{proof}

It is easy to see that any one-dimensional subspace fails Shapiro's Theorem.
Indeed, suppose $(X, \{A_i\})$ is an approximation scheme, and $Y = \spn[e]$ is a
$1$-dimensional subspace of $X$. Let $\varepsilon_n=E(e_1,A_n)$.
Then $\{\varepsilon_n\} \searrow 0$, and every $y=\alpha e\in Y$ satisfies $E(y,A_n)=|\alpha|\varepsilon_n=\mathbf{O}(\varepsilon_n)$. Thus,
Lemma~\ref{sum_of_subsp} implies the following.

\begin{corollary}\label{countable}
Suppose $(X,\{A_i\})$ is an approximation scheme, and $Y$ is a subspace of $X$
with a finite or countable Hamel basis. Then $Y$ fails Shapiro's Theorem relative
to $(A_i)$. In particular:
\begin{enumerate}
\item Any finite dimensional subspace of $X$ fails Shapiro's Theorem.
\item Any separable subspace of $X$ contains a dense subspace, failing Shapiro's Theorem.
\end{enumerate}
\end{corollary}

\begin{proof}[Proof of Theorem~\ref{fin_codim}]
(1)
Suppose $Y$ is a finite codimensional (not necessarily closed) subspace of $X$.
Let $E\subset X$ be a finite dimensional subspace of $X$ such that $Y+E=X$.
By Lemma~\ref{sum_of_subsp} and Corollary \ref{countable}, if $Y$ fails
Shapiro's Theorem with respect to $(A_n)$, then $X=Y+E$ also fails Shapiro's
Theorem with respect to $(A_n)$, which contradicts our assumptions.

(2)
Obviously, the sets  $A_n/Y$, $n=0,1,\ldots$ (which denote
the images of the sets $A_n$ under the quotient map $X\to X/Y$) form an
approximation scheme on $X/Y$.
A direct application of \cite[Proposition 2.1]{almira_oikhberg} shows that
there exists $N\in \N$ such that $A_N + Y = X$.
Consider an approximation scheme in $X$, consisting of sets $B_k = A_{K(N+k-1)}$
($k \in \N$).
By Theorem~\ref{teo_introd}(f), there exists $x \in X$
such that $E(x,B_k) \geq \varepsilon_k$ for $k \geq 1$.
In particular, $E(x, A_{K(N)}) \geq  \varepsilon_1$, and
$E(x, A_{K(n)}) \geq \varepsilon_n$ for any $n > N$.
Find $y \in Y$ such that $x - y = a \in A_N$.
For $n \leq N$, we see that
$E(y,A_n) \geq E(x, A_{K(N)}) > \varepsilon_1\geq \varepsilon_n$, while
for $n > N$, $E(y,A_n) \geq E(x, A_{K(n)}) > \varepsilon_n$.
\end{proof}

\begin{remark}\label{fin_codim_needed}
The assumption of $Y$ being finite codimensional is essential
in Proposition~\ref{fin_codim_p} and Theorem~\ref{fin_codim}.
Indeed, for any infinite codimensional closed subspace $Y$ of
a separable Banach space $X$, there exists an approximation scheme
$A_0 \subset A_1 \subset A_2 \subset \ldots \subset X$, which satisfies
Shapiro's Theorem in $X$, but such that $Y$ fails Shapiro's Theorem
relative to $(A_n)$. To construct such a scheme, recall that a collection
$(e_i)_{i \in I}$ of elements of a Banach space $E$ is called a {\it complete
minimal system} if $\overline{\spn[e_i : i \in I]} = E$, and, for every
$j \in I$, $e_j \notin \overline{\spn[e_i : i \neq j]}$.
Every separable Banach space contains a complete minimal
system (see \cite[Theorem 1.27]{HMVZ} for a stronger result).
Pick a complete minimal system ${\mathcal{D}}$ in $X/Y$. For $n \in \N$,
define $A_n$ as the set of $x \in X$ for which $q(x)$ ($q : X \to X/Y$
denotes the quotient map) can be represented as a linear combination of
no more than $n$ elements of ${\mathcal{D}}$. It follows from
\cite[Theorem 6.2]{almira_oikhberg} that $(X, \{A_n\})$ satisfies Shapiro's Theorem.
However, $Y \subset A_1$.

A more interesting example can be given if $X$ is $\ell_p$ ($0 < p < \infty$)
or $c_0$. Suppose $1 = \varepsilon_0 \geq \varepsilon_1 \geq \ldots \geq 0$.
Then $X$ contains a linear approximation scheme $(A_k)$ and a subspace $Y$,
such that (i) $A_k$ is isometric to $X$ for any $k$, (ii) $Y$ is isometric
to $X$, and (iii) $E(y, A_k) = \varepsilon_k \|y\|$ for any $k \geq 0$,
and any $y \in Y$.

We can view $X$ as a closed linear span of unit vectors $(e_{ij})_{i,j \in \N}$,
with $\|\sum_{ij} \alpha_{ij} e_{ij}\| = (\sum_{ij} |\alpha_{ij}|^p)^{1/p}$.
Set $A_0 = \{0\}$. For $k \geq 1$ define
$A_k = \overline{\spn[e_{ij} : 1 \leq i \leq k , j \in \N]}$, and let
$\gamma_k = (\varepsilon_{k-1}^p - \varepsilon_k^p)^{1/p}$.
For $j \in \N$ set $f_j = \sum_i \gamma_i e_{ij}$, and let
$Y = \overline{\spn[f_j : j \in \N]}$. Then any $y \in Y$ can be
represented as $y = \sum_j \alpha_j f_j = \sum_{ij} \alpha_j \gamma_i e_{ij}$, with
$$
\|y\| = \Big( \sum_{ij} |\alpha_j|^p \gamma_i^p \Big)^{1/p} =
\Big(\sum_j |\alpha_j|^p\Big)^{1/p} ,
$$
hence $Y$ is isometric to $X$. Furthermore, for such $y$,
$$
E(y,A_k) = \Big( \sum_{j=1}^\infty |\alpha_j|^p \sum_{i=k+1}^\infty \gamma_i^p \Big)^{1/p} =
\Big( \sum_{j=1}^\infty |\alpha_j|^p \Big)^{1/p}
\Big( \sum_{i=k+1}^\infty \gamma_i^p \Big)^{1/p} = \varepsilon_k \|y\| .
$$
\end{remark}

Furthermore, the property of satisfying
Shapiro's Theorem is not stable under contractive embeddings.

\begin{proposition}\label{subspace_shapiro}
Suppose $X$ is a separable Banach space.
Then there exists a conti\-nuous embedding of $Z = \ell_1((0,1])$
into $X$, and a family $(A_i)$ in $Z$, such that:
\begin{enumerate}
\item
$(A_i)$ is an approximation scheme in both $Z$ and $X$. Moreover, $(Z,\{A_i\})$ satisfies Shapiro's Theorem.
\item
$A_i$ is dense in $X$ for every $i$ (hence $Z$ is dense in $X$).
Consequently, $(X, \{A_i\})$ fails Shapiro's Theorem.
\end{enumerate}
\end{proposition}

\begin{proof}
%
Let $(x_k)_{k=1}^\infty$ be a
complete minimal system in $X$ such that $\|x_k\| < 1/2^k$ for each $k$. Then
the map $\phi : (0,1] \to X$ given by  $t \mapsto \sum_{k=1}^{\infty} t^k x_k$ is continuous.
Moreover, by \cite{klee}, $\phi$ is injective. By
Theorem 1.56 of \cite{HMVZ}, if $t_1, t_2, \ldots$ are distinct, and $c_1, c_2,\ldots $ are such that $\sum_{j=1}^{\infty} |c_j|$
is finite and $\sum_{j=1}^{\infty} c_j \phi(t_j) = 0$, then $c_j = 0$ for every $j$
(the theorem is stated for basic sequences $(x_k)$, but it works for
minimal systems, too). Finally, also by \cite{klee}, if $(t_j)$ is a sequence convergent to $0$, then $\spn[\{\phi(t_j)\}]$
is dense in $X$.

Take $Z = \ell_1((0,1])$. Denote the ``canonical'' basis in $Z$ by
$(e_t)$ for $0 < t \leq 1$. Define $J : Z \to X$ by setting $J(e_t) = \phi(t)$ and extending it by linearity.
It follows from the properties of $\{x_k\}$ and $\phi$ that $J$ is injective. Moreover,
\[
\|J(e_t)\|_X=\|\phi(t)\|_X\leq \sum_{k=1}^{\infty}\left(\frac{t}{2}\right)^k=\frac{1}{1-t/2}-1=\frac{t/2}{1-t/2}\leq 1 \ \ (t\in (0,1]),
\]
so that $J$ is bounded (hence continuous). Finally, set $A_i$ to be the closed linear span of all $e_t$, for $t$ not in
the set $\{1/i, 1/(i+1), \ldots\}$. Then clearly $(A_i)$ is a nontrivial linear approximation
scheme in $Z$, hence it satisfies Shapiro's Theorem. However, $J(A_i)$ is a
dense linear subspace of $X$.
\end{proof}

Finally, we observe a connection between
Shapiro's Theorem for subspaces,
and some fundamental results of approximation theory. The classical theorem of
Jackson shows that any ``sufficiently smooth'' function is ``well approximable''
(see e.g. \cite[Chapter 7]{devore}). To study this phenomenon in the abstract setting,
suppose $(A_n)$ is an approximation scheme in a quasi-Banach space $X$, and $Y$ is
a quasi-semi-Banach space, continuously and strictly included in $X$.
We say that the approximation scheme $(X,\{A_n\})$ satisfies a
({\it generalized}) {\it Jackson's Inequality} with respect to $Y$ if
there exists a sequence $(c_n)$ such that 
$\lim_{n\to\infty}c_n=0$, and $E(y,A_n)\leq c_n\|x\|_Y$ for all $y\in Y$.
In the classical case of $X = C({\mathbb{T}})$, $A_n=\mathcal{T}_n$
(the set of trigonometric polynomials of degree $\leq n$), and
$Y = C^r({\mathbb{T}})$, we can take $c_n = \gamma_r n^{-r}$.

Suppose $(A_n)$ is an approximation scheme in $X$, and $Y$ is
continuously embedded into $X$. Then $Y$ fails Shapiro's Theorem
relative to $(A_n)$ if and only if there exists a function $C : Y \to [0,\infty)$,
and a sequence $\{\varepsilon_n\}\searrow 0$ such that
$E(y,A_n) \leq \varepsilon_n C(y)$ for all $n$ and $y$.
Thus, the failure of $Y$ to satisfy Shapiro's Theorem relative to $(A_n)$ can be
viewed as a weak form of Jackson's inequality. In fact, we have:

\begin{proposition}\label{far_jackson}
Suppose $(A_n)$ is an approximation scheme in a quasi-Banach space $X$,
and a quasi-Banach space $Y$ is continuously embedded into $X$. Then
the following are equivalent:
\begin{enumerate}
\item
The approximation scheme $(A_n)$ satisfies a
Jackson's inequality with respect to $Y$.
\item
There is no $c > 0$ so that $Y$ is $c$-far from $(A_n)$.
\end{enumerate}
\end{proposition}

\begin{proof}
$(1) \Rightarrow (2)$:
Suppose $\lim_n c_n = 0$, and the inequality $E(y,A_n) \leq c_n \|y\|_Y$
holds for any $y \in Y$ and $n \in \mathbb{N}$. Then
$E(S(Y),A_n) = \sup_{\|y\|=1} E(y,A_n) \leq c_n$, hence
$Y$ cannot be far from $(A_n)$.

\noindent $(2) \Rightarrow (1)$:
Let $c_n = E(S(Y),A_n)$. By assumption, $\lim_n c_n = 0$.
Then, for any $y \in Y$ and $n \in \mathbb{N}$,
$E(y,A_n) = E(y/\|y\|_Y, A_n) \|y\|_Y \leq c_n \|y\|_Y$,
yielding $(1)$.

\end{proof}

\section{Complemented subspaces}\label{complemsub}

Suppose $(A_n)$ is an approximation scheme in a
quasi-Banach space space $X$, and $P$ is a bounded projection
from $X$ onto its subspace $Y$ (clearly, $Y$ is closed).
Then $(Y,\{P(A_n)\})$ is an approximation scheme,
and it is natural to ask under which conditions
$Y$ satisfies Shapiro's Theorem with respect to $(P(A_n))$.
A partial answer is given in Theorem~\ref{subsp_shapiro}.
In particular, we show that, if $Y$
satisfies Shapiro's Theorem with respect to $(P(A_n))$, then it also
satisfies Shapiro's Theorem with respect to $(A_n)$.
If $Y$ is a closed finite codimensional subspace of $X$, and
$P$ is a bounded projection from $X$ onto $Y$, then $(Y,\{P(A_n)\})$
satisfies Shapiro's theorem whenever $(X,\{A_n\})$ does
(Theorem \ref{1_far_for_stable_schemes}).
As an intermediate step for the proof of this last result, we prove that approximation
schemes $(X,\{A_n\})$ satisfying Shapiro's Theorem are stable under the addition of
finite dimensional subspaces of $X$ (Theorem \ref{stability}).

\begin{theorem}\label{subsp_shapiro}
Suppose $P$ is a bounded projection from a quasi-Banach space $X$ onto its
closed subspace $Y$, and $(A_n)$ is an approximation scheme in $X$.
The following are equivalent:
\begin{enumerate}
\item
$Y$ satisfies Shapiro's Theorem with respect to $(P(A_n))$.
\item
There exists a constant $c > 0$ and an infinite set $N_0 \subset \N$
such that, for any $n \in N_0$, there exists
$y \in Y \backslash \overline{P(A_{K(n)})}$
satisfying $E(y, P(A_n))\leq cE(y, P(A_{K(n)}))$.
\item
There is no sequence $\{\varepsilon_n\}\searrow 0$ such that
$E(y,P(A_n))\leq \varepsilon_n\|y\|$ for all $y\in Y$ and $n\in\N$.
\end{enumerate}
Moreover, if $Y$ satisfies Shapiro's Theorem with respect to $\{P(A_n)\}$,
then it also satisfies Shapiro's Theorem with respect to $(A_n)$.
Finally, if $Y$ is Banach and satisfies Shapiro's Theorem with respect to $\{P(A_n)\}$, then for every $\{\varepsilon_n\}\searrow 0$ there exists
an element $y\in Y$ such that $E(y,A_n)\geq \varepsilon_n$ for $n=0,1,2,\ldots$.
\end{theorem}

\begin{proof}
For $n \in \N$, define $B_n=\overline{P(A_n)}$. By assumption,  $\cup_n B_n\supseteq \cup_n P(A_n)$ is dense in $Y$, so $(Y,\{B_n\})$ is an approximation scheme (the other properties of an approximation scheme are inherited from
$(A_n)$). The first part of the theorem  follows from Theorem \ref{teo_introd}, parts (a), (b), and (c) (see also \cite[Theorem 2.2]{almira_oikhberg}). The rest of the theorem follows from part $(f)$ of the same theorem (see also \cite[Corollary 3.6]{almira_oikhberg}) and the fact that for any $a \in A_k$,
$\|P\|\|y-a\|\geq \|P(y-a)\|=\|y-Pa\|\geq E(y,B_k)$.
\end{proof}
In general, an infinite dimensional subspace of $X$ needs not satisfy
Shapiro's Theorem (see Remark~\ref{fin_codim_needed}). However,
certain subspaces do satisfy it.

\begin{corollary}\label{complemented}
Suppose $P$ is a bounded projection from a quasi-Banach space $X$ onto its
closed subspace $Y$. Suppose, furthermore, that $(A_n)$ is a non-trivial linear approximation scheme on $X$
(i.e., $K(n)=n$ and $\overline{A_n}\neq X$ for all $n\in\N$)
and $Y\not\subseteq \bigcup_{n\in\N}\overline{P(A_n)}$.
Then $Y$ satisfies Shapiro's Theorem relative to $(A_n)$. If, in addition,
$Y$ is a Banach space, then for any sequence $\{\varepsilon_n\}\searrow 0$ there exists
$y \in Y$ such that $E(y,A_n) \geq \varepsilon_n$ for any $n$.
\end{corollary}
\begin{proof}
The condition $Y\not\subseteq \bigcup_{n\in\N}\overline{P(A_n)}$ guarantees that $(Y,\{P(A_n)\})$ is a non-trivial linear approximation scheme, so that it satisfies Shapiro's Theorem.
\end{proof}

\begin{remark}
Corollary \ref{complemented} is not true for arbitrary (non-linear)
approximation schemes $(A_n)$ in $X$ such that $(X,\{A_n\})$ satisfies Shapiro's Theorem.
To see this, consider the following example. Let $(Z,\{Z_n\})$ be an approximation scheme
that satisfies Shapiro's Theorem.  Let $(Y,\{Y_n\})$ be an approximation scheme that fails
Shapiro's Theorem and such that $Y\not\subseteq \bigcup_n\overline{Y_n}$
(there are examples of this in \cite[Section 4]{almira_oikhberg}), let $X=Z\oplus Y$ with quasi-norm $\|(z,y)\|_X=\max\{\|z\|_Z,\|y\|_Y\}$ (hence $P:X\to Y$ given by $P(z,y)=y$ is our projection, $\|P\|=1$). Our approximation scheme is $(X,\{A_n\})$, where $A_n=Z_n+Y_n$. It is clear that this approximation scheme satisfies Shapiro's Theorem, that  $Y\not\subseteq \bigcup_n\overline{P(A_n)}$ and $Y$ fails Shapiro's Theorem  with respect to $\{A_n\}$.
\end{remark}

\begin{remark}\label{projections}
It may happen
that $Y$ satisfies Shapiro's Theorem relative to a linear approximation scheme
$(A_n)$ in the ambient space $X$, but not relative to $(P(A_n))$ ($P$ is a projection
from $X$ onto $Y$). Indeed, consider a Hilbert space $X$ with an orthonormal
basis $e_1, f_1, e_2, f_2, \ldots$. Let $Y = \spn[e_1, e_2, \ldots]$. For $k \geq 1$
define $g_k = k^{-1} e_k + \sqrt{1 + k^{-2}} f_k$, and set
$A_n = \spn[e_1, f_1, \ldots, e_n, f_n, g_{n+1}, g_{n+2}, \ldots]$.
Then $E(e_m, A_n) = \sqrt{1 - m^{-2}}$ for $m > n$, hence, by Theorem~\ref{far_shap},
$Y$ satisfies Shapiro's Theorem relative to $(A_n)$.
On the other hand, $P(A_n)$ is dense in $Y$ for every $n$.
\end{remark}

We next show that the property of satisfying Shapiro's Theorem is stable under
adding a finite dimensional subspace.

\begin{theorem}\label{stability}
Suppose an approximation scheme $(X,\{A_n\})$ satisfies Shapiro's Theorem, and
$F$ is a finite dimensional subspace of $X$. Then the approximation scheme $(X,\{A_n+F\})$
also satisfies Shapiro's Theorem.
\end{theorem}

\begin{proof}
We assume, with no loss of generality, that $X$ is $p$-convex. Suppose, for
the sake of contradiction, that $(X,\{A_n+F\})$ fails Shapiro's Theorem.
Assume first that $F\cap (\bigcup \overline{A_n})=\{0\}$.
By Theorem \ref{teo_introd}, the fact that $(X,\{A_n+F\})$ fails Shapiro's Theorem
implies the existence $N_0\in\N$ such that $E(S(X),A_{N_0}+F)<\frac{1}{4}$.
Hence for every $x\in S(X)$ there exists $e(x)\in F$ and $a(x)\in A_{N_0}$ such that $\|x-a(x)+e(x)\|<2E(S(X),A_{N_0}+F)<\frac{1}{2}$. Thus, $\|a(x)-e(x)\|^p\leq \|a(x)-e(x)-x\|^p+\|x\|^p\leq \frac{2^p+1}{2^p}$. By the finite dimensionality of $F$,
and the fact that $F \cap \overline{A_{N_0}} = \{0\}$, for every $e\in F$ we have
\[
E(e,A_{N_0})=E(\frac{e}{\|e\|},A_{N_0})\|e\|\geq \rho \|e\|,
\]
where $\rho=\inf_{x\in S(F)}E(x,A_{N_0})>0$.  Hence
\[
 \|e(x)\| \leq \frac{1}{\rho} E(e(x),A_{N_0})\leq \frac{1}{\rho}\|e(x)-a(x)\|\leq \frac{1}{\rho}\left(\frac{2^p+1}{2^p}\right)^{\frac{1}{p}}=C<\infty.
\]
The approximation scheme $(X,\{A_n\})$ satisfies Shapiro's Theorem, hence
$E(S(X),A_n)=1$ for all $n\in\N$. Hence, for all $n\geq 1$ we can take $x_n\in S(X)$ such that $E(x_n,A_n)\geq 1-\frac{1}{n}$. The boundedness of the associated sequence $\{e(x_n)\}$ in conjunction with the finite dimensionality of $F$ imply that there exists $e_*\in F$ and a subsequence $e(x_{n_k})$ such that $\|e_*-e(x_{n_k})\|\to 0$ for $k\to\infty$. Take $\varepsilon>0$ such that $2\varepsilon^p<1-1/2^p$. For $k$ big enough we get
\[
\|x_{n_k}-e_*+a(x_{n_k})\|^p\leq \|x_{n_k}-e(x_{n_k})+a(x_{n_k})\|^p+\|e(x_{n_k})-e_*\|^p\leq \frac{1}{2^p}+ \varepsilon^p.
\]
On the other hand, the density of $\bigcup A_n$ implies that there exists $N_1\in \N$ and
$b\in A_{N_1}$ such that $\|b-e_*\|<\varepsilon$. Pick $k > K(\max\{N_0,N_1\})$
so large that $(1-1/n_k)^p > 2^{-p} + 2\varepsilon^p$. Then
\begin{eqnarray*}
\Big(1-\frac{1}{n_k}\Big)^p &\leq & E(x_{n_k},A_{n_k})^p\leq E(x_{n_k},A_{K(\max\{N_0,N_1\})})^p\\
&\leq&  \|x_{n_k}-b+a(x_{n_k})\|^p \leq \|x_{n_k}-e_*+a(x_{n_k})\|^p+\|e_*-b\|^p\leq \frac{1}{2^p}+ 2\varepsilon^p,
\end{eqnarray*}
yielding a contradiction.

In the general case, note that $\bigcup \overline{A_n}$ is a linear subspace
of $X$, hence $F_0 = F\cap (\bigcup \overline{A_n})$ is a subspace of $F$.
One can see that there exists $n_0 \in \N$ such that
$F_0 = F\cap \overline{A_n}$ for $n \geq n_0$.
Find a subspace $F_1$ of $F$ such that $F_1 \cap F_0 = \{0\}$, and
$F_1 + F_0 = F$. Then $F_1 \cap \bigcup \overline{A_n}=\{0\}$, and
$A_{K(n)} + F_1 \supset A_n + F$ for $n \geq n_0$. The family $B_n = A_{K(n)}$
forms an approximation scheme in $X$. By Theorem~\ref{teo_introd},
$(B_n)$ satisfies Shapiro's Theorem whenever $(A_k)$ does.
By the reasoning above,
the approximation scheme $(X,\{A_{K(n)}+F_1\})$ satisfies Shapiro's Theorem.
Therefore, so does the original approximation scheme $(X,\{A_n+F\})$.
\end{proof}


Recall that any finite codimensional closed subspace $Y$ of a quasi-Banach space $X$ is complemented.
Indeed, there exists a finite dimensional subspace $F$ of $X$, such that $F \cap Y = \{0\}$,
and $X = Y + F$. Any $x \in X$ has a unique representation $x = y + f$, with $y \in Y$
and $f \in F$. We can define a projection $Q$ from $X$ onto $F$ by setting
$Q(x) = f$. It is easy to see that $Q$ is bounded, hence so is $P = I - Q$.
It follows that $P$ is a bounded projection from $X$ onto $Y$.

\begin{theorem} \label{1_far_for_stable_schemes}
Suppose $(X,\{A_n\})$ satisfies Shapiro's Theorem, and $P$ is a bounded
projection onto a closed finite codimensional subspace $Y$ of $X$. Then $Y$
satisfies Shapiro's Theorem with respect to $(P(A_n))$.
Moreover, 
$E(S(Y),A_n)\geq \frac{1}{\|P\|}$. Consequently,
if $X$ is a Hilbert space and $Y$ is a finite codimensional
closed subspace of $X$, then $Y$ is $1$-far from the approximation scheme $(A_n)$.
\end{theorem}

This result provides an improvement over Proposition~\ref{fin_codim_p} when $\|P\|$ is small.
Note that the existence of a projection $P$ as above follows from the paragraph
preceding the proposition.

\begin{proof}
Recall that there exists a constant $C_q > 0$ such that $\|x_1 + x_2\| \leq C_q(\|x_1\| + \|x_2\|)$
for any $x_1, x_2 \in X$.

Let $Q = I - P$, and $F = Q(X)$.
Then $\bigcup_nP(A_n)$ is dense in $Y$, so that $(Y,\{P(A_n)\})$ is an approximation scheme.
Note that $P(A_n) + Q(A_n) \subset B_n = A_n + F$ for every $n$.
Indeed, fix $a, b \in A_n$. Then $Pa +Qb = a + Q(b-a)$, and $Q(b-a) \in F$.

Assume, for the sake of contradiction, that $(Y,\{P(A_n)\})$ does not satisfy Shapiro's Theorem.
Then there exists $\{\varepsilon_n\}\in c_0$ such that for every $y\in Y$ and $n \geq 0$,
$E(y,P(A_n))\leq \varepsilon_n C(y)$. For any $x \in X$, we have
$$
E(x,B_n) \leq E(Px + Qx, P(A_n) + Q(A_n)) \leq C_q(E(Px, P(A_n)) + E(Qx, F)) .
$$
However, $E(Qx, F) = 0$, hence $E(x,B_n) \leq \varepsilon_n C_q C(Px)$ for every $x$.
The desired contradiction arises when we recall that, by Theorem~\ref{stability},
$(X, \{B_n\})$ satisfies Shapiro's Theorem.

Consequently, $E(S(Y),P(A_n))=1$ for all $n\in \N$. This, in conjunction with the inequality
$\|P\|\|y-a\|\geq \|P(y-a)\|=\|y-Pa\|\geq E(y,P(A_k))$ (which holds
for $y\in Y$ and $a\in A_k$), implies that, for any $y\in Y$,
$E(y,A_k)\geq\frac{1}{\|P\|}E(y,P(A_k))$, so $E(S(Y),A_k)\geq \frac{1}{\|P\|}$.
If $X$ is a Hilbert space, then $\overline{Y}$ is $1$-far from $(A_n)$.
The density of $Y$ inside $\overline{Y}$ implies that $Y$ is also $1$-far
from $(A_n)$.
\end{proof}

\section{Boundedly compact approximation schemes}\label{compact_schemes}

This section is devoted to the approximation schemes $(A_n)$ which are
{\it boundedly compact} in $X$ -- that is, the set
$\overline{\{a \in A_n : \|a\| \leq 1\}}$ is compact for every $n$.
In this case, infinite dimensional closed subspaces satisfy
Shapiro's Theorem (Theorem~\ref{compact_shap}). For some schemes $(A_n)$,
an even stronger statement holds (Theorem~\ref{bern_subspace}).
These results are then used to study approximability of analytic
functions (Corollary~\ref{suavidad}).

To proceed, we need an auxiliary result.

\begin{lemma}\label{compact_far}
Suppose $(A_i)$ is a boundedly compact approximation scheme in
a $p$-convex quasi-Banach space $X$.
Then any infinite dimensional subspace of $X$ is $2^{-1/p}$-far from $(A_i)$.
If, moreover, $X$ is a Banach space (that is, $X$ is $1$-convex), any
infinite dimensional subspace of $X$ is $1$-far from $(A_i)$.
\end{lemma}

An application of Theorem~\ref{far_shap} yields:

\begin{theorem}\label{compact_shap}
Suppose $Y$ is an infinite dimensional  closed subspace of a quasi-Banach space $X$, and
the approximation scheme $(A_n)$ is boundedly compact in $X$.
Then $Y$ satisfies Shapiro's Theorem relative to $(A_n)$.
\end{theorem}

\begin{proof}[Proof of Lemma~\ref{compact_far}]
Consider first the case of $X$ being a Banach space.
Suppose, for the sake of contradiction, that an infinite dimensional $Y \subset X$
is not $1$-far from $(A_i)$ (we can assume that $Y$ is closed).
Then there exists $n \in \N$ and $c \in (0,1)$ such that for every
$y \in B(Y)$ there exists $a \in A_n$ such that
$\|y - a\| < c$. By the triangle inequality, $\|a\| \leq 2$.

Pick $d \in (c,1)$. By compactness, there exists a finite $(d-c)$-net
$(a_i)_{i=1}^N \subset \{a \in A_n : \|a\| \leq 2\}$.
For any $y \in B(Y)$, there exists $i$ such that  $\|y - a_i\| \leq d < 1$.
Letting $E = \spn[a_i : 1 \leq i \leq N]$, we see that $dist(y, E) \leq d \|y\|$
for any $y \in Y$. This, however, is impossible, by \cite[Lemma 1.19]{HMVZ}.

Now suppose $X$ is quasi-Banach. Suppose, for the sake of contradiction,
there exists $n \in \N$ and $c \in (0,2^{-1/p})$ such that for every
$y \in B(Y)$ there exists $a \in A_n$ such that $\|y - a\| < c$.
Pick $d \in (c,2^{-1/p})$ and $\delta > 0$ satisfying $\delta^p + c^p < d^p$.
Suppose $(a_i)$ is a $\delta$-net in $\{a \in A_i : \|a\| \leq 2^{1/p}\}$.
We claim that, for every $y \in B(Y)$, there exists $\ell$ such that
$\|y - a_\ell\| < d$. Indeed, pick $a \in A_i$ such that
$\|y - a\| < c$. By $p$-convexity, $\|a\|^p \leq \|y\|^p + \|y - a\|^p < 2$.
Find $\ell$ to satisfy $\|a - a_\ell\| < \delta$. By our choice
of $\delta$, $\|y - a_\ell\| < d$.

Note that, for every infinite-dimensional quasi-Banach space $Z$, and every
$\lambda < 1$, there exists a sequence $(z_i)_{i \in \N}$ such that
$\|z_i - z_j\| > \lambda$ whenever $i \neq j$. Indeed, it is well
known (see e.g.~\cite[Lemma 6.3]{almira_oikhberg}) that, if $E$
is a proper closed subspace of a quasi-Banach space $F$, then there
exists $f \in B(F)$ such that $dist(f,E) > \lambda$. We use this
fact to construct $(z_i)$ inductively: pick an arbitrary norm $1$ $z_1$.
If $z_1, \ldots, z_k$ with the desired properties have already been constructed,
find $z_{k+1} \in B(Z)$ such that $dist(z_{k+1}, \spn[z_1, \ldots, z_k]) > \lambda$.

Thus, there exists a sequence $(y_i)_{i \in \N}$ such that
$\|y_i - y_j\| > 2^{1/p} d$ whenever $i \neq j$. However, there exist distinct
$i$ and $j$ such that $\|y_i - a_\ell\| < d$ and $\|y_j - a_\ell\| < d$,
for some $\ell$. Then $\|y_i - y_j\|^p \leq \|y_i - a_\ell\|^p + \|y_j - a_\ell\|^p < 2 d^p$,
which is a contradiction.
\end{proof}

More can be said when the approximation scheme in question is linear (it is
easy to see that a linear approximation scheme $(A_n)$ is boundedly
compact if and only if $\dim A_n < \infty$ for every $n$).

\begin{theorem}\label{bern_subspace}
Suppose $\{0\} = A_0 \subset A_1 \subset A_2 \ldots$ is a sequence of finite
dimensional subspaces of a Banach space $X$, $Y$ is an infinite dimensional closed subspace
of $X$, and $\{\varepsilon_n\} \searrow 0$. Then there exists $y \in Y$
such that $\|y\| = \varepsilon_0$, and $E(y,A_n)_X \geq \varepsilon_n$ for
any $n \geq 0$.
\end{theorem}

This is a generalization of the classical Bernstein's Lethargy Theorem.
Note that, in general, we cannot guarantee the existence of $y \in Y$
with the property that $E(y, A_n) = \varepsilon_n$. For instance, suppose
$X = \ell_2$ (with the canonical basis $e_1, e_2, \ldots$),
$A_n = \spn[e_1, \ldots, e_n]$, and $Y = \spn[e_3, e_4, \ldots]$.
Then $E(y, A_1) = E(y,A_2)$ for any $y \in Y$. Moreover, the hypothesis of
$Y$ being a closed subset of $X$ can not be deleted in the theorem.
Indeed, $Y=\bigcup_nA_n$ is an infinite dimensional subspace of $X$,
and for every $y\in Y$, $E(y,A_n)=0$ for sufficiently large $n$.

\begin{proof}
We briefly sketch the proof, using the ideas of \cite[pp.~264-266]{singerlibro}.
Inductively, we can construct a sequence of finite dimensional subspaces
$0 = \{0\} \subset B_1 \subset B_2 \subset \ldots \subset X$ in such a way
that, for every $k$, $A_k \subset B_k$, and $B_k \cap Y \not\subset B_{k-1}$
(here we use the fact that $\dim Y = \infty$). Then we construct, for each
$n \geq 0$, $y_n \in B_{n+1} \cap Y$ satisfying $E(y_n, B_k) = \varepsilon_k$ for
$0 \leq k \leq n$. To this end, fix $n$, and find
$z_n \in B_{n+1} \cap Y$ for which $E(z_n, B_n) = \varepsilon_n$.
This can be done, since $B_{n+1} \cap Y \not\subset B_{n}$, so that there exists $z\in B_{n+1} \cap Y$ with $E(z,B_n)>0$, and now it is easy to find $\lambda>0$ such that $\phi(\lambda)=E(\lambda z,B_n)=|\lambda|E(z,B_n)=\varepsilon_n$. Take $z_n=\lambda z$. Pick $w_n \in (B_n \cap Y) \backslash B_{n-1}$. Then there exists
$\lambda_n \in \R$ such that  $E(z_n + \lambda_n w_n, B_{n-1}) = \varepsilon_{n-1}$.
Set $z_{n-1} = z_n + \lambda_n w_n$. Note that, as $w_n \in B_n$,
$E(z_{n-1}, B_n) = E(z_n, B_n) = \varepsilon_n$.
On the next step, we obtain $z_{n-2} = z_{n-1} + \lambda_{n-1} w_{n-1}$,
for some $w_{n-1} \in (B_{n-1} \cap Y) \backslash B_{n-2}$ and $\lambda_[n-1] \in \R$,
such that $E(z_{n-2}, B_n) = \varepsilon_n$, $E(z_{n-2}, B_{n-1}) = \varepsilon_{n-1}$,
and $E(z_{n-2}, B_{n-2}) = \varepsilon_{n-2}$. Proceeding further in the same manner,
we end up with $z_0 \in B_{n+1} \cap Y$, satisfying $E(z_0, B_k) = \varepsilon_k$
for $0 \leq k \leq n$ (in particular, $\|z_0\| = \varepsilon_0$). Let $y_n = z_0$.

For $0 \leq k \leq n$, pick $u_{nk} \in B_k$ satisfying
$\|y_n - u_{nk}\| = \varepsilon_k$. Clearly $\|u_{nk}\| \leq 2 \varepsilon_0$.
Using compactness and diagonalizing (as on p.~265 of \cite{singerlibro}),
find $n_1 < n_2 < \ldots$ such that  the sequence $(u_{n_j k})_{j=1}^\infty$
converges for every $k$. We claim that the sequence $(y_{n_j})$ converges
to $y \in Y$, satisfying $E(y, B_k)_X = \varepsilon_k$ for every $k$.
It suffices to show that, for every $\delta > 0$, there exists
$N \in \N$ such that  $\|y_{n_i} - y_{n_j}\| < \delta$ whenever $i, j > N$.
To this end, pick $k$ so large that $\varepsilon_k < \delta/3$. Pick
$N \in \N$ such that  $\|u_{n_i k} - u_{n_j k}\| < \delta/3$ for any
$i, j > N$. By the triangle inequality, $\|y_{n_i} - y_{n_j}\| < \delta$
for such $i$ and $j$.
\end{proof}

\begin{remark}
It is important to note that Theorem \ref{bern_subspace} does not follow from Corollary \ref{complemented}, since there are examples of infinite dimensional closed subspaces $Y$ of a Banach $X$ such that there is no bounded projection $P:X\to Y$ (we say that $Y$ is uncomplemented in $X$).
It is  well known  that every closed subspace $Y$ of $X$ which is finite dimensional or finite codimensional, is complemented. In 1971 J. Lindenstrauss
and L. Tzafriri \cite{LZ} proved that if every closed subspace of a Banach space $X$ is complemented, then $X$ is isomorphic to a Hilbert space.
A classical example of an uncomplemented subspace is provided by
$X=C(\T)$, $Y=A(\D)$ (the disk algebra -- see \cite{hoffman}).  Another elementary example is $X=\ell_{\infty}$ and $Y=c_0$.
T. Gowers and B. Maurey \cite{GM} constructed a Banach space $X$ such that every closed subspace $Y$ of $X$  which is not finite dimensional nor finite codimensional, is uncomplemented in $X$. 
\end{remark}


We apply Theorems \ref{compact_shap} and \ref{bern_subspace} to the study
of real analytic functions on an interval.
if we consider $C^r[a,b]$ ($r \geq 1$) as a subspace of $C[a,b]$.
Let $A_n$ be the space of algebraic polynomials of degree not exceeding $n$.
A classical theorem of Jackson (see e.g.~\cite[Theorem 8.6.2]{devore})
shows that, for $f \in C^r[a,b]$, $E(f,A_n) = \mathbf{O}(n^{-r})$.
Below, we show the speed of decay of the sequence $(E(f,A_n))$ can no longer
be controlled if the smoothness of $f$ is violated at $a$ and $b$, but
$f$ is analytic on $(a,b)$. That is, the conditions
guaranteeing that a function is ``well approximable'' must be
``of global nature'' (holding on the whole domain).

\begin{corollary}\label{suavidad}
Suppose $(A_n)$ is an approximation scheme on $C[a,b]$. Then
\begin{itemize}
\item[$(a)$] If $(A_n)$ is boundedly compact, then for all
$\{\varepsilon_n\}\searrow 0$ there exists  $f\in C[a,b]$ which is analytic
in $(a,b)$, such that $E(f,A_n)\not=\mathbf{O}(\varepsilon_n)$.
\item[$(b)$] If the sets $A_n$ are finite dimensional subspaces of $C[a,b]$,
then for all $\{\varepsilon_n\}\searrow 0$ there exists $f\in C[a,b]$, analytic in
$(a,b)$, such that $E(f,A_n)\geq \varepsilon_n$ for $n=0,1,2,\ldots$.
\end{itemize}
\end{corollary}

\begin{proof}
Given $M\in\mathbb{R}$, the operator $T_M:C[a,b]\to C[a+M,b+M]$ given by
$T(f)(x)=f(x-M)$ is a linear isometry of Banach spaces. In particular, $T_M$
preserves relatively compact sets and finite dimensional subspaces. Moreover,
$(A_n)$ is an approximation scheme on $C[a,b]$ if and only if $(T_M(A_n))$
is an approximation scheme on $C[a+M,b+M]$. Finally, $f\in C[a,b]$ is real analytic
at $\alpha\in (a,b)$ if and only if $T_M(f)$ is real analytic at $\beta=\alpha+M$.

Thus, we can assume that $0 < a < b$.
Consider $Y=\overline{\mathbf{span}[\{x^{n^2}:n\in\N\}]}^{C[a,b]}$.
By M\"{u}ntz Theorem (see \cite[Theorem 11]{almira}), $Y$ is a proper
subspace of $C[a,b]$. Furthermore, by Full Clarkson-Erd\"{o}s-Schwartz Theorem
(see \cite[Theorems 28 and 31]{almira}), the elements of $Y$ have analytic extensions
to the set $\{z\in\C\setminus (-\infty,0]: a<|z|<b\}$.
An application of Theorem~\ref{compact_shap} (or
Theorem~\ref{bern_subspace}) establishes (a) (respectively,~(b)).
\end{proof}

\section{Compactly embedded subspaces}\label{sec_compact}

In this section we investigate the case when $Y$ is compactly embedded into
$X$ (that is, the unit ball of $Y$ is relatively compact in $X$).
Theorem~\ref{compact_subspaces} shows that, in this case, $Y$ cannot be far
from an approximation scheme $(A_n)$. Consequently, $Y$ fails Shapiro's
Theorem, and moreover, it satisfies Jackson's inequality
(see Proposition \ref{far_jackson}).
Furthermore, by Theorem~\ref{compact_embed}, if $(A_n)$ is boundedly
compact, then the subspaces $Y$ failing Shapiro's Theorem with respect to $(A_n)$
are precisely those satisfying $Y\subseteq Z$ for a
certain space $Z$, compactly embedded into $X$.
We also provide examples of compactly embedded subspaces.

\begin{theorem}[Jackson's theorem for compact embeddings]
\label{compact_subspaces}
Suppose $(A_n)$ is an approximation scheme in $X$, and $Y$ is a subspace of $X$,
such that the inclusion $Y\hookrightarrow  X$ is compact. Then
$Y$ is not far from $(A_n)$. Consequently, $Y$ fails Shapiro's Theorem
relative to $(A_n)$.
\end{theorem}

\begin{proof}
Assume, for the sake of contradiction, that $Y$ is far from $(A_n)$.
That is, $$\inf_{n\in\N}E(S(Y),A_n)=2c>0.$$ Then there exists a sequence
$\{y_n\}_{n=0}^{\infty} \subset S(Y)$ such that $E(y_n,A_n)>c$ for all $n\in\N$.
Now, the compactness of the inclusion $Y\hookrightarrow X$, when applied to
the sequence $\{y_{K(n)}\}_{n=0}^{\infty}$, implies that there exists a sequence
$n_i\to\infty$ and an element $y\in X$ such that
$\lim_{i\to\infty}\|y_{K(n_i)}-y\|=0$. Hence
\[
E(y_{K(n_i)},A_{K(n_i)})\leq C_q[E(y_{K(n_i)}-y,A_{n_i})+E(y,A_{n_i})]\to 0,
\]
which contradicts the fact that $c<E(y_{K(n_i)},A_{K(n_i)})$ for all $i$.
The failure of Shapiro's Theorem then follows by Theorem~\ref{far_shap}(1).
\end{proof}

\begin{remark} Let $Y$ be a subspace of $X$, and let $W\subseteq Y$ be a
homogeneous subset of $Y$ (i.e., $\lambda W\subseteq W$ for all scalars $\lambda$).
If $S(Y)\cap W$ is a relatively compact subset of $X$, then the same arguments of
Theorem \ref{compact_subspaces} (changing $Y$ by $W$ and $S(Y)$ by $W\cap S(Y)$)
prove that, for each approximation scheme $(A_n)$ in $X$ there exists a
sequence $\{\varepsilon_n\}\in c_0$ (depending on $(A_n)$) such that for all $w\in W$,
$E(w,A_n)=\mathbf{O}(\varepsilon_n)$. \end{remark}

To illustrate the scope of Theorem~\ref{compact_subspaces}, we provide a few
examples of compactly embedded subspaces. The first one is a simple application
of Ascoli-Arzela Theorem.

\begin{example} Let $(A_n)$ be any approximation scheme on $C[a,b]$, and $Y$ is either
$C^{(1)}[a,b]$ or $\mathbf{Lip}_{\alpha}[a,b]$ ($\alpha>0$). Then
$Y$ is compactly embedded into $C[a,b]$.
\end{example}

Now consider the space $BV(\Omega)$ of functions of bounded variation
on $\Omega$. To be more precise,
suppose $\Omega$ is an open subset of $\mathbb{R}^N$. Let
$$BV(\Omega)=\{u\in L^1(\Omega): \|u\|_{BV(\Omega)}:=
 \sup_{\phi \in \mathbf{C}_c^{(\infty)}(\Omega,\mathbb{R}),
 \sup_{x\in\Omega}|\phi(x)|\leq 1}\int_{\Omega} u(x)\mathbf{div}(\phi)(x)dx<\infty\}$$
be the space of functions of bounded variation on $\Omega$.
Equipping $BV(\Omega)$ with the norm $\|u\|=\|u\|_{L^1(\Omega)}+\|u\|_{BV(\Omega)}$,
we turn it into a Banach space.
Furthermore, the embedding $BV(\Omega)\hookrightarrow L^1(\Omega)$ is compact
(see e.g.~\cite[Chapter 3]{ambrosio}).

\begin{example}\label{cor_convex_L1}
Let $(A_n)$ be any approximation scheme in $L^1(\Omega)$.
Then there exists a sequence $\{\varepsilon_n\}\searrow 0$ such that
$E(f,A_n)_{L^1(\Omega)}={\mathbf{O}}(\varepsilon_n)$
for any $f\in BV(\Omega)$.
Consequently, if $(A_n)$ is an approximation scheme in $L^1(a,b)$, then
there exists a sequence $\{\varepsilon_n\}\searrow 0$ such that
$E(f,A_n)_{L^1(\Omega)}={\mathbf{O}}(\varepsilon_n)$ whenever $f$
is a bounded monotone or convex function on $(a,b)$.
\end{example}

\begin{proof}
As noted above, the embedding of $BV(\Omega)$ into $L^1(\Omega)$ is compact.
The general result now follows from Theorem~\ref{compact_subspaces}.
In the particular case of $\Omega = (a,b)$, it is well known that any
bounded monotone function has bounded variation (and conversely, any function
of bounded variation is a  difference of two bounded monotone functions).
Furthermore, for any convex function on $(a,b)$ there exists a $c \in (a,b)$
such that the restrictions of $f$ to $(a,c)$ and $(c,b)$ are monotone, hence
convex functions must have bounded variation.
\end{proof}

\begin{example}\label{integrals}
Let $(A_n)$ be any approximation scheme in $L^1(a,b)$. Define $B_n=\{f(t)=\int_a^tg(s)ds:g\in A_n\}$.
Then $(B_n)$ is an approximation scheme in $C_0[a,b]=\{f\in C[a,b]:f(a)=0\}$ and there exists a sequence $\{\varepsilon_n\}\searrow 0$ such that, for any convex function
$f\in C_0[a,b]$, $E(f,B_n)_{C[a,b]}={\mathbf{O}}(\varepsilon_n)$.
\end{example}

\begin{proof}
To show that $(B_n)$ is an approximation scheme in $C_0[a,b]$,
it suffices to show that $\cup_n B_n$ is dense in $C_0[a,b]$.
It is easy to see that polynomials vanishing at $a$ are dense in $C_0[a,b]$,
hence it suffices to show that, for any such polynomial $p$, and any $\varepsilon > 0$,
there exists $f \in B_n$ with $\|p - f\| < \varepsilon$. To this end, find
$g \in A_n$ such that $\|p^\prime - g\|_{L^1} < \varepsilon/(b-a)$.
Then the function $f(t) = \int_a^t g(s) \, ds$ ($a \leq t \leq b$) belongs to $B_n$.
Furthermore, for $a \leq t \leq b$, $p(t) = \int_a^t p^\prime(s) \, ds$,
hence
$$
\|p - f\| \leq \sup_t \int_a^t |p^\prime(s) - g(s)| \, ds < \varepsilon .
$$

To prove the second statement, take into account that if
$f\in C_0[a,b]$ is convex, then $f(t)=\int_a^tg(s)ds$ for a certain increasing function $g\in L^1(a,b)$. Hence we can use Example~\ref{cor_convex_L1} to prove that there exists a sequence $\{\varepsilon_n\}\searrow 0$ such that
\begin{eqnarray*}
E(f,B_n)_{C_0[a,b]} &=& \inf_{a_n\in A_n}|f(t)-\int_a^ta_n(s)ds|\\
&=&\inf_{a_n\in A_n}|\int_a^tg(s)ds-\int_a^ta_n(s)ds| \leq \inf_{a_n\in A_n}\int_a^t|g(s)-a_n(s)|ds\\
&\leq& E(g,A_n)_{L^1(a,b)}=\mathbf{O}(\varepsilon_n).
\end{eqnarray*}
\end{proof}

Still another application of Theorem~\ref{compact_subspaces} is in order:

\begin{example}\label{fourier}
Suppose $\{\phi_k\}_{k=0}^{\infty}$ is an orthonormal basis in a separable
Hilbert space $H$. For $x\in H$ and $k\in\mathbb{N}$, denote by
$c_k(x)= \langle x,\phi_k \rangle$ the $k$-th Fourier coefficient of
$x$ with respect to $\{\phi_k\}_{k=0}^{\infty}$. Let $\{c_k^*(x)\}$ stand for
the non-increasing rearrangement of $\{|c_k(x)|\}$. Let $Y\subset H$ be a subspace of $H$ which is compactly embedded into $H$. Then there exists a decreasing sequence $\{\varepsilon_n\}_{n=0}^{\infty}\searrow 0$ such that, for all $y\in Y$,
\[
c_n^*(y) \leq \Big( \sum_{k=n}^{\infty}c_k^*(y)^2\Big)^{\frac{1}{2}}   =\mathbf{O}(\varepsilon_n).
\]
\end{example}

\begin{proof}
Let $A_n=\bigcup_{\{i_1,i_2,\ldots,i_n\}\subseteq \mathbb{N}}\mathbf{span}[\{\phi_{i_k}\}_{k=1}^n]$, $n=1,2,\ldots$.
Then $(H,(A_n))$ is an approximation scheme, and
$$
\Big( \sum_{k=n}^{\infty}c_k^*(y)^2 \Big)^{\frac{1}{2}}=E(y,A_{n-1}) .
$$
We complete the proof by applying Theorem~\ref{compact_subspaces}.
\end{proof}

\begin{theorem}\label{compact_embed}
Suppose $(X,\{A_n\})$ is a boundedly compact approximation scheme, and
$Y$ is continuously embedded subspace of $X$.
Then:
\begin{enumerate}
\item
$Y$ fails Shapiro's Theorem with respect to $(A_n)$ if and only if
there exists a quasi-Banach space $Z\subseteq X$ such that $Y\subseteq Z$,
and the embedding $Z \hookrightarrow X$ is compact.
\item
If $Y$ is not far from $(A_n)$,
then $Y$ is compactly embedded into $X$.
\end{enumerate}
\end{theorem}

\begin{proof}
(1)
The property of failing Shapiro's Theorem is inherited by subspaces.
If $Z$ is compactly embedded into $X$, it fails Shapiro's Theorem
by Theorem \ref{compact_subspaces}. In this situation, $Y$ will also
fail Shapiro's Theorem.

Conversely, if $Y$ fails Shapiro's Theorem with respect to $(A_n)$,
then there exists a sequence $\{\varepsilon_n\}\in c_0$, such that $E(y,A_n)=\mathbf{O}(\varepsilon_n)$ for every $y \in Y$. By
\cite[Lemma 2.3]{almira_oikhberg}, we may assume that
$\varepsilon_n\leq 2\varepsilon_{K(n+1)-1}$ for all $n\in \N$.
Then ${\mathbf{A}}(\varepsilon_n)=\{x\in X:
 \|\{\frac{E(x,A_n)}{\varepsilon_n}\}\|_{\ell_{\infty}}<\infty\}$
is a quasi-Banach subspace of $X$
(see \cite[Remark 3.5 and Proposition 3.8]{almiraluther2}).
Moreover, ${\mathbf{A}}(\varepsilon_n)$ satisfies the generalized
Jackson's inequality
 $E(y,A_n)\leq \varepsilon_n\|y\|_{{\mathbf{A}}(\varepsilon_n)}$.
By Proposition~\ref{far_jackson}, the space ${\mathbf{A}}(\varepsilon_n)$
cannot be far from $(A_n)$. By part
$(2)$ of this theorem (see also \cite[Theorem 3.32]{almiraluther2}),
the natural inclusion of ${\mathbf{A}}(\varepsilon_n)$ to $X$ is compact.
To complete the proof, take $Z={\mathbf{A}}(\varepsilon_n)$.

(2)
We show that, for every $c > 0$, $X$ contains a
finite $c$-net for $B(Y) = \{y \in Y : \|y\|_Y \leq 1\}$.
Without loss of generality, we can assume
that $Y$ is embedded into $X$ contractively. Recall that there exists
a constant $C_q$ so that $\|x_1 + x_2\|_X \leq C_q(\|x_1\|_X + \|x_2\|_X)$
for any $x_1, x_2 \in X$. Let $\varepsilon_n = E(S(Y),A_n)$.
By assumption, $\lim_n \varepsilon_n = 0$.
Pick $n$ so large that $\varepsilon_n < c/(2C_q)$.
By compactness, there exists a finite $c/(2C_q)$-net $\{a_1, \ldots, a_N\}$
in $\{a \in A_n : \|a\|_X \leq 2C_q\}$. We show that $\{a_1, \ldots, a_N\}$
is also a $c$-net for $B(Y)$. Indeed, for any $y \in B(Y)$, there
exists $a \in A_n$ such that $\|y - a\|_X \leq c/(2C_q)$.
As $\|a\|_X \leq C_q(\|y\|_X + \|y - a\|_X) < 2 C_q$,
hence there exists $\ell \in \{1, \ldots, N\}$ such that
$\|a - a_\ell\|_X < c/(2C_q)$. Then $\|y - a_\ell\|_X \leq
C_q (\|y - a\|_X + \|a - a_\ell\|_X) \leq c$, and we are done.
\end{proof}

\begin{remark}\label{completeness_needed}
If the assumptions of Theorem \ref{compact_embed}$(2)$ are satisfied,
then, by Theorem \ref{compact_subspaces}, $Y$ fails Shapiro's Theorem.
Conversely, if $Y$ is a closed subspace of $X$, failing Shapiro's Theorem
relative to $(A_n)$, then, by Theorem \ref{far_shap}, $Y$ can not be far
from $(A_n)$. In this situation, by Theorem \ref{compact_embed}$(2)$,
$Y$ is compactly embedded. However, non-closed subspaces $Y$ of $X$,
failing Shapiro's Theorem relative to a boundedly compact scheme $(A_n)$,
need not be compactly embedded. As an example, consider
the space $Y$, described in Remark~\ref{part_3_theorem21},
equipped with the norm inherited from $X = C[0,2\pi]$.
The sets ${\mathcal{A}}_n$, consisting of all algebraic
polynomials of degree less than $n$, form a
boundedly compact approximation scheme in $X$.
$Y$ contains $\cup_n {\mathcal{A}}_n$, hence it is dense in $X$, and
its embedding into $X$ is not compact. Furthermore, $Y$ is a
proper subspace of $X$, hence it is not complete.
By Remark~\ref{part_3_theorem21}, $Y$ fails Shapiro's Theorem
relative to $({\mathcal{A}}_n)$.
\end{remark}


\section{$1$-far subspaces of finite codimension}\label{1-far}

Suppose $(A_i)$ is an approximation scheme in $X$. By Theorem~\ref{far_shap},
a closed subspace $Y \subset X$ satisfies Shapiro's Theorem relative to
$(A_i)$ is $c$-far from $(A_i)$, for some $c \in (0,1]$. In this section, we
investigate the extremal case of $1$-far subspaces. Recall that, by
Theorem~\ref{teo_introd}, if $(X,\{A_i\})$ satisfies Shapiro's Theorem,
then $X$ is $1$-far from $(A_i)$. The main result of this
section is Theorem \ref{basis_far}, describing a class spaces $X$,
so that every finite codimensional $Y \subset X$ is $1$-far from
an approximation scheme $(A_i)$, provided $(A_i)$ satisfies Shapiro's Theorem.
In particular, $X = (\sum_{i \in I} \ell_{p_i})_{p_0}$ is such a space,
provided $1 < \inf_{i \in I \cup \{0\}} p_i \leq
\sup_{i \in I \cup \{0\}} p_i < \infty$ (Corollary~\ref{l_p}).
The main tool in our investigation is the Defining Subspace Property (DSP),
introduced in Definition~\ref{def:dsp}. This property may be of further interest
to Banach space experts, and is studied throughout this section.

Recall that the {\it modulus of convexity} of a Banach space $X$ is defined by setting,
for $0 < \varepsilon \leq 2$,
\begin{equation}
\moco_X(\varepsilon) = \inf \Big\{ 1 - \frac{\|x+y\|}{2} :
\|x\| \leq 1, \, \|y\| \leq 1 , \|x-y\| \geq \varepsilon \Big\} .
\label{mod_conv}
\end{equation}
Clearly, the function $\moco_X$ is non-decreasing.
A Banach space is called {\it uniformly convex} if $\moco_X(\varepsilon) > 0$
for any $\varepsilon \in (0,2)$. It is known (see e.g. \cite[Section 1.e]{LT2})
that $L^p$ spaces are uniformly convex for $1 < p < \infty$.
Moreover, any uniformly convex space is reflexive \cite[Proposition 1.e.3]{LT2}.

A Banach space $X$ is said to have the {\it Reverse Metric Approximation Property}
({\it RMAP} for short) if, for any finite dimensional subspace $F$ of $X$, and for
any $\delta > 0$, there exists a finite rank operator $u \in \mathbf{B}(X)$ such that
$u|_F = I_F$, and $\|I_X - u\| < 1 + \delta$. $X$ has the {\it shrinking RMAP} if, for
any finite dimensional subspaces $F \subset X$ and $G \subset X^*$, and any
$\delta > 0$, there exists a finite rank $u \in \mathbf{B}(X)$ satisfying
$u|_F = I_F$, $\|u^*|_G - I_G\| < \delta$, and $\|I_X - u\| < \delta$.
By a small perturbation argument, in both definitions above we can only
require $\|u|_F - I_F\| < \delta$.
The reader is referred to \cite{CK} for more information about the RMAP.

We also need to introduce a new definition, reflecting the mutual position
of finite codimensional and finite dimensional spaces.

\begin{definition}\label{def:dsp}
Suppose $Y$ is a closed finite codimensional subspace of a Banach space $X$,
$F$ is a finite dimensional subspace of $X$, and $\delta > 0$.
We say that $F$ is {\it $(\varepsilon,\delta)$-defining} for $Y$ if any $x \in X$
with $\|x\| \leq 1$ and $E(x,F) > 1 - \delta$, we have $E(x,Y) \leq \varepsilon$.
$Y$ has the {\it Defining Subspace Property} ({\it DSP} for short)
if for every $\varepsilon > 0$ there exist $\delta > 0$, and a
$(\varepsilon,\delta)$-defining finite dimensional subspace.
\end{definition}

The DSP can be thought of as a generalization of orthogonality.
Indeed, suppose $Y$ is an infinite dimensional subspace of a Hilbert space $X$.
Let $F = Y^\perp$. Then, for any $x \in X$ with $\|x\| \leq 1$,
$E(x,Y)^2 = \|x\|^2 - E(x,F)^2 \leq 1 - E(x,F)^2$.
Thus, $F$ is $(\varepsilon,\sqrt{1-\varepsilon^2})$-defining for $Y$,
for any $\varepsilon > 0$.

Let us now state the main result of this section.
\begin{theorem}\label{basis_far}
The following statements hold:
\begin{itemize}
\item[$(1)$] Suppose $X$ is a Banach space, and an approximation scheme $(X, \{A_i\})$
satisfies Shapiro's Theorem. Suppose, furthermore, that $Y$ is a
finite codimensional subspace of $X$, with the Defining Subspace Property.
Then $Y$ is $1$-far from $(A_i)$.
\item[$(2)$] Suppose $X$ is a uniformly convex Banach space with the Reverse
Metric Approximation Property. Then any
finite codimensional subspace of $X$ has the Defining Subspace Property.
\end{itemize}
Consequently, if $X$ is a uniformly convex Banach space, with the Reverse
Metric Approximation Property and the approximation scheme $(X, \{A_i\})$ satisfies Shapiro's Theorem,
then any finite codimensional subspace of $X$ is $1$-far from $(A_i)$.
\end{theorem}

\begin{remark}\label{rem:fin_codim}
Note that, in Theorem \ref{basis_far}, we do not make any assumptions about
the nature of $(A_i)$, only about the geometry of $X$. For particular
schemes $(A_i)$, $Y$ can be shown to be $1$-far from $(A_i)$.
For instance, by \cite[Lemma 6.4]{almira_oikhberg}, any  finite codimensional
subspace of a Banach space is $1$-far from a linear approximation
scheme $(A_i)$ if $\dim X/\overline{A_i} = \infty$ for any $i$.
More examples of $1$-far subspaces are given in Lemma~\ref{compact_far}, and
Theorems~\ref{1_far_for_stable_schemes}, \ref{gen_haar}, \ref{uniform}, and
\ref{uncond_far}.

Note also that, by Theorem \ref{teo_introd}, $X$ is $1$-far from $(A_i)$
whenever $X$ satisfies Shapiro's Theorem relative to $(A_i)$.
Furthermore, by Proposition~\ref{fin_codim_p}, any finite codimensional subspace of $X$
is $1/2$-far from $(A_i)$. We do not know whether such a subspace must
be $1$-far from $(A_i)$.
\end{remark}

Theorem~\ref{basis_far} implies:

\begin{corollary}\label{l_p}
Consider an index set $\Gamma$, and sets $({\mathcal{F}}_i)_{i \in \Gamma}$
such that $\cup_{i \in \Gamma} {\mathcal{F}}_i$ is infinite.
Suppose the family $(p_i)$ {\rm{(}}$i \in \Gamma^\prime = \Gamma \cup \{0\}${\rm{)}}
satisfies
$1 < \inf_{i \in \Gamma^\prime} p_i \leq \sup_{i \in \Gamma^\prime} p_i < \infty$.
Suppose, furthermore, that an approximation scheme $(A_i)$ in
$X = \big(\sum_{i \in \Gamma} \ell_{p_i}({\mathcal{F}}_i) \big)_{p_0}$
satisfies Shapiro's Theorem. Then any finite codimensional subspace
of $X$ is $1$-far from $(A_i)$.
\end{corollary}

For the proof, recall that a Banach lattice $X$ is called {\it $p$-convex}
(resp. {\it $p$-concave}) if there exists a constant $C$ such that the
inequality $\|(\sum |x_i|^p)^{1/p}\| \leq C (\sum \|x_i\|^p)^{1/p}$
(resp. $(\sum \|x_i\|^p)^{1/p} \leq C \|(\sum |x_i|^p)^{1/p}\|$)
holds for any collection $x_1, \ldots, x_n \in X$. The infimum of all
$C$'s for which the above inequalities hold is denoted by $M^{(p)}(X)$
(resp. $M_{(p)}(X)$), and is called the {\it $p$-convexity}
(resp. {\it concavity}) constant of $X$. The reader is referred to
\cite[Section 1.d]{LT2} for more information on these notions.
To give just one example, an application of Minkowski Inequality shows
that the Banach lattice $L^r$ is $u$-convex and $v$-concave, with
constant $1$, whenever $1 \leq u \leq r \leq v \leq \infty$.

\begin{proof}[Proof of Corollary~\ref{l_p}]
To show that $X$ has the RMAP, consider the set
${\mathcal{F}} = \{ (i, \alpha) : i \in \Gamma, \alpha \in {\mathcal{F}}_i\}$.
Then, for $x = (x_{i \alpha})_{(i, \alpha) \in {\mathcal{F}}}$,
$$
\|x\|^{p_0} = \sum_i \big( \sum_\alpha |x_{i \alpha}|^{p_i} \big)^{p_0/p_i} .
$$
If $F$ is a finite subset of ${\mathcal{F}}$, define a projection $P_F$ by setting
$(P_F x) = x_{i \alpha}$ if $(i,\alpha) \in F$,
$(P_F x)x_{i \alpha} = 0$ otherwise. Clearly, $I_X - P_F$ is contractive, hence
$X$ has the RMAP.

To prove the uniform convexity of $X$,
let $p = \min\{2, \inf_{i \in \Gamma^\prime} p_i\}$ and
$q = \max\{2, \sup_{i \in \Gamma^\prime} p_i\}$. As noted in the paragraph preceding
this proof, $\ell_{p_i}({\mathcal{F}}_i)$  is $p$-convex and $q$-concave with
constant $1$. Therefore, $M^{(p)}(X) = M_{(q)}(X) = 1$. By \cite[Theorem 1.f.1]{LT2},
$X$ is uniformly convex. To complete the proof, apply Theorem~\ref{basis_far}.
\end{proof}

As we shall see below, general $L^p$ spaces may fail the RMAP.

\begin{proof}[Proof of Theorem~\ref{basis_far}(1)]
Let $n\in\mathbb{N}$.  By hypothesis, given $\varepsilon>0$ there exists a finite dimensional space $F$ and $\delta \in (0, \varepsilon)$ such that
$E(x,F)>1-\delta$ and $\|x\|\leq 1$ imply $E(x,Y)<\varepsilon$. On the other hand, Theorem~\ref{stability} implies that there exists $x_*\in S(X)$ such that $\min\{E(x_*,F),E(x_*,A_n)\} \geq E(x_*,A_n+F)\geq 1-\delta$. Hence  $E(x_*,Y)<\varepsilon$. Take $y_*\in Y$ such that $\|x_*-y_*\|<2\varepsilon$. Then $\|y_*\|\leq 1+2\varepsilon$ and
\[
E(y_*,A_n)\geq E(x_*,A_n)-\|x_*-y_*\| \geq 1-\delta-2\varepsilon\geq 1-3\varepsilon.
\]
This shows that $E(S(Y),A_n)=1$ since $\varepsilon>0$ was arbitrary.
\end{proof}

To prove part (2) of Theorem~\ref{basis_far}, we need an auxiliary result.



\begin{lemma}\label{reflexive}
Any reflexive Banach space with the RMAP has the shrinking RMAP.
Consequently, a reflexive Banach space has the RMAP if and only
if its dual has it.
\end{lemma}

\begin{proof}
Suppose $F$ and $G$ are finite dimensional subspaces of $X$ and $X^*$,
respectively, and $\delta \in (0,1)$. By the definition of the RMAP, there
exists a net $(u_\alpha)_{\alpha \in {\mathcal{A}}}$ of finite
rank operators on $X$, such that $u_\alpha|F = I_F$,
$u_\alpha \to I_X$ pointwise, and $\|u_\alpha\| < 1 + \delta$.
Then $u_\alpha^* \to I_{X^*}$ in the point-weak$^*$ topology.
As $X$ is reflexive, we conclude that $u_\alpha^* x^* \to x^*$
weakly, for any $x^* \in X^*$.

Pick a $\delta/9$-net $g_1, \ldots, g_N$ in the unit ball of $G$. Consider
$\tilde{g} = (g_1, \ldots, g_N) \in \tilde{X} = \ell_\infty^N(X^*)$,
and the maps $\tilde{u}_\alpha$ on $\tilde{X} = \ell_\infty^N(X^*)$, taking
$(x_i^*)_{i=1}^N$ to $(u_\alpha^* x_i^*)_{i=1}^N$.
Then $\tilde{u}_\alpha \tilde{x} \to \tilde{x}$ in the weak topology of $\tilde{X}$.
In particular, $\tilde{x}$ belongs to the weak closure of the set
$\{\tilde{u}_\alpha \tilde{x}\}_{\alpha \in {\mathcal{A}}}$.
Applying Mazur' Theorem \cite[Appendix F]{AK} to the set
$\{\tilde{u}_\alpha \tilde{g}\}_{\alpha \in {\mathcal{A}}}$,
we find $\alpha_1, \ldots, \alpha_m$ and $\lambda_1, \ldots, \lambda_m \in (0,1)$,
such that $\sum_k \lambda_k = 1$, and
$\|\sum_k \lambda_k \tilde{u}_{\alpha_k} \tilde{g} - \tilde{g}\| < \delta/2$.
We claim that the operator $u = \sum_k \lambda_k u_{\alpha_k}$ has the desired properties.
Indeed, $u_{|F}=I_F$, and
$$
\|I - u\| \leq \sum_k \lambda_k \|I - u_{\alpha_k}\| < 1 + \delta .
$$
It remains to show that $\|u^* g - g\| \leq \delta \|g\|$ for any $g$ in the
unit ball of $G$. Find $i$ with $\|g_i - g\| < \delta/9$. Then
$$
\|u^* g_i - g_i\| \leq \|\sum_k \lambda_k \tilde{u}_{\alpha_k} \tilde{g} - \tilde{g}\| <
\delta/2 .
$$
Furthermore, $\|u\| < 3$, hence
$$
\|u^* g - g\| \leq \|u^* g_i - g_i\| + \|u^*\| \|g_i - g\| + \|g - g_i\| <
\delta/2 + \delta/3 + \delta/9 < \delta .
$$
\end{proof}

\begin{proof}[Proof of Theorem~\ref{basis_far}(2)]
The space $X$ is uniformly convex, hence reflexive (and even superreflexive).
For $\varepsilon \in (0,1/2)$, pick $\delta \in (0,\varepsilon/2)$ satisfying
$(1-\delta)/(1+\delta) > 1 - \moco_X(\varepsilon/2)$.
Suppose $Y$ is a finite codimensional subspace of $X$.
By Lemma~\ref{reflexive}, $X^*$ has the RMAP, and therefore, there exists a
finite rank $u \in \mathbf{B}(X)$ satisfying $u^*|_{Y^\perp} = I_{Y^\perp}$, and
$\|I_X - u\| < 1 + \delta$. Note that, in this situation, $(I_X - u)(X) \subset Y$.
Indeed, for any $x \in X$ and $x^* \in Y^\perp$,
$\langle x^*, x \rangle = \langle u^* x^*, x \rangle = \langle x^*, u x \rangle$,
hence
$\langle x^*, (I_X - u) x \rangle = \langle x^*, x \rangle - \langle x^*, u x \rangle = 0$.

Let $F = u(X)$, and suppose $1-\delta < E(x,F) \leq \|x\| \leq 1$. It suffices to show that
$\|ux\| \leq \varepsilon$. To this end,
let $y = (I - u)x = x - ux$, and $z = (x+y)/2 = x - ux/2$. As $ux \in F$, we have
$\|z\| \geq E(x,F) > 1 - \delta$. On the other hand, $x^\prime = x/(1+\delta)$ and
$y^\prime = y/(1+\delta)$ belong to the unit ball of $X$, hence
$$
\| \frac{z}{1+\delta} \| = \|\frac{x^\prime + y^\prime}{2}\| \leq
1 - \moco_X(\|x^\prime - y^\prime\|) \leq 1 - \moco_X(\|ux\|/2) .
$$
Thus, $(1-\delta)/(1+\delta) \leq 1 - \moco_X(\|ux\|/2)$, which implies
$\|ux\| \leq \varepsilon$.
\end{proof}

Our next result shows that, in Theorem~\ref{basis_far}$(2)$,
neither the uniform convexity nor the RMAP can be omitted. Moreover, the comments
below Proposition~\ref{no_DSP} show that to satisfy RMAP or DSP is, in general,
a strong geometric assumption.

To proceed further, recall the definition of generalized Schatten spaces.
Suppose $\se$ is a symmetric sequence space. That is, suppose $\se$ is a
Banach space of sequences, such that the sequences with finitely many
non-zero entries are dense in $\se$, and
$\|(x_i)_{i \in \N}\|_\se = \|(\omega_i x_{\pi(i)})_{i \in \N}\|_\se$ whenever
$|\omega_i| = 1$ for every $i$, and $\pi : \N \to \N$ is a bijection.
We say that $\se$ has the {\it Fatou property} if, for any sequence
$x = (x_i)_{i \in \N}$, if
$\sup_n \|(x_1, \ldots, x_n, 0, \ldots)\|_\se < \infty$, then
$x \in \se$, and
$\|x\|_\se = \sup_n \|(x_1, \ldots, x_n, 0, \ldots)\|_\se$.

If $\se$ is a symmetric sequence space with the Fatou property, we
define the Schatten space $\ce$ as the set of those compact operators
$T \in \textbf{B}(\ell_2)$ such that $(s_i(T))\in\se$
(here, $s_1(T) \geq s_2(T) \geq \ldots \geq 0$ are the singular numbers of $T$).
By e.g..~\cite{GK, Si}, $\ce$ becomes a Banach space when endowed with the norm
$\|T\|_{\se} = \|(s_i(T))\|_{\se}$. Furthermore, by \cite[Lemma III.6.1]{GK},
for $T \in \ce$,
\begin{equation}
\inf_{\rank u < n} \|T - u\|_\se = \|(s_n(T), s_{n+1}(T), \ldots)\|_\se .
\label{best_appr}
\end{equation}
Thus, finite rank operators are dense in $\ce$ whenever $\se$ is separable.
Observe that ${\mathcal{S}}_{c_0}$ is just the space $K(\ell_2)$ of
compact operators, while ${\mathcal{S}}_{\ell_p} = \is_p$
($1 \leq p < \infty$) is the usual Schatten $p$ space. The reader is
referred to \cite{Si} or \cite{GK} for more information.

\begin{proposition}\label{no_DSP}
The following Banach spaces have subspaces of codimension $1$, failing the
Defining Subspace Property:
\begin{enumerate}
\item
$c_0$.
\item
$L^p(0,1)$, with $p \in (1,2) \cup (2,\infty)$.
\item
The Schatten space $\ce$, where $\se$ is a symmetric sequence space,
not isomorphic to the Hilbert space, and satisfying
$M^{(p)}(\se) = M_{(q)}(\se) = 1$ for some $1 < p \leq q < \infty$.
\end{enumerate}
\end{proposition}

Clearly, $c_0$ has the RMAP, but it is not uniformly convex.
On the other hand, $L^p$ is uniformly convex for $1< p < \infty$.
By Proposition \ref{no_DSP} and Theorem \ref{basis_far}(2), it fails the RMAP for
$p\neq 2$.
Furthermore, by \cite{TJ84}, $\ce$ is uniformly convex for $\se$ as in (3).
Consequently, $\ce$ fails the RMAP. In fact, a stronger statement is true:
if $\se$ is a reflexive symmetric sequence space, such that $\ce$ embeds isometrically into
a space with the RMAP, then $\ce$ is $4$-isomorphic to a Hilbert space.
Furthermore, a separable rearrangement invariant Banach space
of functions on $(0,1)$ or $(0,\infty)$ embeds isometrically into a space with the RMAP
if and only if it is isometric to a Hilbert space. Both of these facts have been
established in \cite{Oi}.

For the proof of Proposition~\ref{no_DSP}, we need  a few  technical results.
The first one deals with functions on $(0,1)$ and involves
nothing but computations.

\begin{lemma}\label{int_0}
Suppose $p \in[1,\infty] \backslash \{2\}$, and $\alpha \in (0,1/2)$.
Denote by $\one$ the function equal to $1$ everywhere on $(0,1)$.
Then there exist positive numbers $a$ and $b$, and a real number $c$,
so that $\alpha a = (1-\alpha) b$,
$\|a \chi_{(0,\alpha)} - b \chi_{(\alpha,1)}\|_p = 1$, and
$\|a \chi_{(0,\alpha)} - b \chi_{(\alpha,1)} + c \one\|_p < 1$.
\end{lemma}

Next we handle $\ce$. Note that, if $M^{(p)}(\se) = M_{(q)}(\se) = 1$
for some $1 < p \leq q < \infty$, then $\se$ is {\it regular} in the
terminology of \cite[Section 1.7]{Si} -- that is, for any
$(x_i)_{i \in \N} \in \se$, we have
$\lim_n \|(x_n, x_{n+1}, \ldots)\|_\se = 0$.
Denote by $E_{ij}$ the $(i,j)$ matrix unit --
that is, the matrix with $1$ in the $(i,j)$ position, and zeroes elsewhere.
We can identify the dual of $\ce$ with $\is_{\se^\prime}$
(see e.g. \cite[Chapter 3]{Si}) via the parallel duality:
$\langle (b_{ij}), (a_{ij}) \rangle = \sum a_{ij} b_{ij}$.
Then $E_{ij}^* \in \is_{\se^\prime}$ defines a contractive linear functional:
$\langle E^*_{ij}, (a_{uv}) \rangle = a_{ij}$. The projection $P_{ij}$
``onto the $(i,j)$ entry'' can be defined by setting
$P_{ij} a = \langle E^*_{ij}, a \rangle E_{ij}$.
That is, for $a = (a_{uv})$, $(P_{ij} a)_{k\ell} = a_{ij}$ if $k=i$ and $\ell=j$,
and $(P_{ij} a)_{k\ell} = 0$ otherwise.

\begin{lemma}\label{1-1}
Suppose $\se$ is a symmetric sequence space with
$M^{(p)}(\se) = M_{(q)}(\se) = 1$ for some $1 < p \leq q < \infty$, not
isomorphic to $\ell_2$. Then $\|I - P_{11}\|_{\mathbf{B}(\ce)} > 1$.
\end{lemma}

Note that $\|I - P_{11}\|_{\mathbf{B}(\ce)} = \|I - P_{ij}\|_{\mathbf{B}(\ce)}$
for any pair $(i,j)$.

\begin{proof}
Suppose, for the sake of contradiction, that $\|I - P_{11}\|_{\mathbf{B}(\ce)}  = 1$.
Then $\prod_{k=1}^m (I - P_{i_k j_k})$ is contractive for every finite family
$((i_k, j_k))_{k=1}^m$. Therefore, for any $A \subset \N \times \N$,
the projection $Q_A$ is contractive. Here, $Q_A$ is defined by setting
$Q_A E_{ij} = \left\{ \begin{array}{ll}
   E_{ij}   &    (i,j) \in A    \cr
   0        &    (i,j) \notin A
\end{array} \right.$. This, in turn, implies that $E_{ij}$ is an unconditional
basis for $\ce$. By \cite[Theorem 2.2]{KP} (or \cite{Pis}), this is only
possible if $\se$ is isomorphic to a Hilbert space.
\end{proof}

The following simple lemma deals with small perturbations of finite dimensional subspaces.

\begin{lemma}\label{small_pert}
Suppose $\delta \in (0,1)$, and $F$ and $F^\prime$ are finite dimensional
subspaces of a Banach space $X$, such that $E(S(F),F^\prime) < \delta$.
Then, for every $x \in X$ with $\|x\| \leq 1$,
$E(x,F^\prime) \leq E(x,F) + 2\delta$.
\end{lemma}

\begin{proof}
For any $c > 0$, there exists $f \in F$ with $\|x-f\| < E(x,F) + c$.
By the triangle inequality, $\|f\| < 2 + c$. Find $f^\prime \in F^\prime$
with $\|f - f^\prime\| < (2+c)\delta$. Then
$$
E(x,F^\prime) \leq \|x-f^\prime\| \leq \|x-f\| + \|f - f^\prime\| <
E(x,F) + c + (2+c)\delta .
$$
We conclude the proof by noting that $c$ can be arbitrarily small.
\end{proof}


We also need a simple observation.

\begin{lemma}\label{dense_sub}
Suppose a finite codimensional subspace $Y$ of a Banach space $X$ has the
Defining Subspace Property. Suppose, furthermore, that
$(F_\alpha)_{\alpha \in {\mathcal{A}}}$ is a family of finite dimensional
subspaces of $X$, such that for any finite dimensional subspace $F$ of $X$,
$\inf_\alpha E(S(F),F_\alpha) = 0$. Then, for every $\varepsilon > 0$,
there exists $\alpha \in {\mathcal{A}}$ and $\delta > 0$ such that
$F_\alpha$ is $(\varepsilon,\delta)$-defining for $Y$.
\end{lemma}

\begin{proof}
There exists finite dimensional $F \subset X$ and $\delta > 0$ such that
$F$ is $(\varepsilon,3\delta)$-defining for $Y$ -- that is, $E(x,Y) \leq \varepsilon$
whenever $1 - 3 \delta < E(x,F) \leq \|x\| \leq 1$. Find $\alpha$ such that
$E(S(F),F_\alpha) < \delta$, and show that $F_\alpha$ is
$(\varepsilon,\delta)$-defining for $Y$.

Suppose $1 - \delta < E(x,F_\alpha) \leq \|x\| \leq 1$. By Lemma~\ref{small_pert},
$1 - 3 \delta < E(x,F)$, hence $E(x,Y) \leq \varepsilon$.
\end{proof}

\begin{proof}[Proof of Proposition~\ref{no_DSP}]
(1)
Denote by $P_k$ the $k$-th basis projection (corresponding to the canonical
basis on $c_0$), and let $Y = \{x \in c_0 : P_1 x = 0\}$. We show that,
for any finite dimensional subspace $F$ of $c_0$, there exists
$x \in c_0$ such that $E(x,Y) = 1 = \|x\|$, and $E(x,F) \geq c$.
Indeed, pick $\delta \in (0,(1-c)/2)$, and find $m \geq n = \dim F$ for which
$\|(P_m - I)|_F\| < \delta$. We claim that $x = (1, \ldots, 1, 0, 0, \ldots)$
($m+1$ $1$'s) has the desired property. The equality $E(x,Y) = 1 = \|x\|$ is
clearly satisfied. Suppose, for the sake of contradiction, that $E(x,F) < c$.
Then there exists $f \in F$ with $\|x-f\| < c$. By the triangle inequality,
$\|f\| \leq \|x\| + \|x-f\| < 2$. Then
$$
\|x - P_m f\| \leq \|x-f\| + \|(P_m - I)|_F\| \|f\| < 1 ,
$$
which contradicts the fact that $E(x, P_m(x)) = 1$.

(2) Let $Y$ be the set of all $g \in L^p(0,1)$ with $\int g = 0$.
By Lemma~\ref{int_0}, there exists $\kappa \neq 0$, $\alpha \in (0,1/2)$,
and positive $a$ and $b$, for which the function
$f = a \chi_{(0,\alpha)} - b \chi_{(\alpha,1)}$ is such that $\int f = 0$,
$\|f\| > 1$, and $1 = \|f - \kappa \one\| = \inf_\gamma \|f - \gamma \one\|$
(in other words, $\alpha |a-\kappa|^p + (1-\alpha) |b+\kappa|^p \leq
\alpha |a-\gamma|^p + (1-\alpha) |b+\gamma|^p$ for any $\gamma$).
Suppose, for the sake of contradiction, that $Y$ has the DSP.
By Lemma~\ref{dense_sub}, we can assume that there exist
$0 = t_0 < t_1 < \ldots < t_n = 1$, such that
$F = \spn[\chi_{(t_{k-1}, t_k)}] : 1 \leq k \leq n]$ is
$(|\kappa|/2,\delta)$-defining for $Y$, for some $\delta$.
Define the function
$$
x = \sum_{k=1}^n \big( a \chi_{(t_{k-1}, \alpha t_{k-1} + (1-\alpha) t_k)} -
b \chi_{(\alpha t_{k-1} + (1-\alpha) t_k, t_k)} \big) - \kappa \one =
\sum_{k=1}^n x_k ,
$$
where
$$
x_k = a \chi_{(t_{k-1}, \alpha t_{k-1} + (1-\alpha) t_k)} -
b \chi_{(\alpha t_{k-1} + (1-\alpha) t_k, t_k)} - \kappa \chi_{(t_{k-1}, t_k)} .
$$
Then $\|x\| = (\sum_{k=1}^n \|x_k\|^p)^{1/p} = 1$. Furthermore, $x_k$ is
supported on $(t_{k-1}, t_k)$. By our choice of $\kappa$,
$$
E(x_k, \spn[\chi_{(t_{k-1}, t_k)}])^p =
\inf_c \|x_k - c \chi_{(t_{k-1}, t_k)}\|^p = \|x_k\|^p = t_k - t_{k-1} .
$$
Furthermore, $E(x,F)^p = \sum_{k=1}^n E(x_k, \spn[\chi_{(t_{k-1}, t_k)}])^p = 1$.
One the other hand, $\int$ is a contractive functional,
vanishing on $Y$. Therefore, $E(x,Y) \geq |\int x| = |\kappa|$, which yields
a contradiction.

(3) We show that $Y = (I - P_{11})(\ce)$ fails the DSP. To this end,
find a norm one matrix $a = (a_{ij})$ with $\kappa = a_{11} > 0$, and such that
$\|a + \gamma E_{11}\| \geq 1$ for any $\gamma$, and $\|a - a_{11} E_{11}\| > 1$
(this is possible, by Lemma~\ref{1-1}). By the discussion preceding
Lemma~\ref{1-1}, $\se$ is separable, hence finite rank operators are
dense in $\ce$. Therefore, matrices with finitely many non-zero entries
are dense in $\ce$. If $Y$ has the DSP, then,
by Lemma~\ref{dense_sub}, there exists $n \in \N$ and $\delta > 0$ such that
$1 - \delta < E(x,F) \leq \|x\| \leq 1$ implies $E(x,Y) \leq \kappa/2$,
with $F = \spn[E_{ij} : 1 \leq i, j \leq n]$. To obtain a contradiction,
consider
$$
x = a_{11} + \sum_{i>1} a_{i1} E_{i+n,1} + \sum_{j>1} a_{1j} E_{1,j+n} +
\sum_{i,j>1} a_{ij} E_{i+n,j+n} .
$$
Then $\|x\|=1$. To show that $E(x,F) = 1$, define the projection $Q$
by setting $Q E_{ij} = E_{ij}$ if either $i \in \{1, n+1, n+2, \ldots\}$
or $j \in \{1, n+1, n+2, \ldots\}$, and $Q E_{ij} = 0$ otherwise. Then
$Q$ is contractive, $Qx=x$, and $QF = \spn[E_{11}]$. Therefore,
$$
E(x,F) = \inf_{f \in F} \|x-f\| = \inf_{f \in F} \|x-Qf\| =
\inf_\gamma \|x - \gamma E_{11}\| = 1 .
$$
Finally, the contractive functional $E_{11}^*$ vanishes on $Y$, hence
$E(x,Y) \geq \langle E_{11}^*, x \rangle = \kappa$, a contradiction.
\end{proof}

We conclude this section by noting that the DSP is ``very fragile.''

\begin{proposition}\label{destroy_DSP}
Suppose $X_0$ is a subspace of an infinite dimensional Banach space $X$
of codimension $1$. Then $X$ can be equivalently renormed in such a way that
$X_0$ has no Defining Subspace Property.
\end{proposition}

\begin{proof}
The space $Y = Z \oplus_\infty X_0$, where $Z$ is a $1$-dimensional space,
can clearly be regarded as a renorming of $X$. We show that, for any finite
dimensional subspace $F$ of $Y$, there exists $y \in Y$ such that
$E(y,F)_Y = 1 = \|y\|$, and $E(y,X_0)_Y = 1$. Enlarging $F$ if necessary,
we can assume that $F = Z \oplus F_0$, where $F_0$ is a finite dimensional
subspace of $X_0$. By \cite[Lemma 1.19]{HMVZ}, there exists $w \in X_0$, for which
$E(w,F_0) = 1 = \|w\|$. It is easy to see that $y = 1 \oplus w$ has all the desired properties.
\end{proof}

\section{Additional examples}\label{examples}

In this section we investigate specific approximation schemes $(A_i)$, so that,
for a wide class of subspaces $Y$ of the ambient space $X$, Shapiro's Theorem is
satisfied. In working with systems of functions, we need the notion of a
generalized Haar family.


\begin{definition}\label{def_haar}
Let $\Omega$ be a topological space and let, for each $n$, $A_n$ be a set of continuous complex valued functions on $\Omega$.
We say that the family $\{A_n\}$ is {\it generalized Haar}
if there exists a function $\psi = \psi_{\{A_n\}} : \N \to \N$ such that
no non-zero function of the form $\Re g$ ($g \in A_n$) has more than
$\psi(n) - 1$ zeroes on $\Omega$.
An approximation scheme $(X,\{A_n\})$ is called {\it generalized Haar}
if  $\{A_n\}$ is a generalized Haar system.
\end{definition}

\begin{definition}\label{def_dictionary}
Suppose ${\mathcal{D}}$ is a subset of a quasi-Banach space $X$.
For $n \in \N$, we define
$$
\Sigma_n(\mathcal{D}) = \cup_{F \subset {\mathcal{D}}, |F| \leq n} \spn[F] .
$$
${\mathcal{D}}$ is called {\it dictionary} if $\overline{\mathbf{span}[\mathcal{D}]}=X$.
${\mathcal{D}}$ is called a {\it generalized Haar system} if the  family $\{\Sigma_n(\mathcal{D})\}$ is generalized Haar.
\end{definition}

These concepts generalize the concept of being a Haar system, which appears when we impose $\psi(n)=n$ for all $n$.

Below we call a subspace $Y$ of $C_0(I)$ {\it pseudo-real} if, for any $f \in Y$,
$\Re f$
belongs to $Y$. In the real case, any
subspace of $C_0(I)$ is pseudo-real. In the complex case, $Y$ is pseudo-real if and only if
$Y = Y_r + i Y_r$, where $Y_r = \{\Re f : f \in Y\}$. In particular, the
span of a family of real-valued functions is pseudo-real.

\begin{theorem}\label{gen_haar}
Suppose $\{A_n\}$ is a generalized Haar system in $C_0(I)$ ($I$ is an interval),
and $Y$ is an infinite dimensional pseudo-real subspace of $C_0(I)$.
Then $Y$ {\rm{(}}equipped with the norm of $C_0(I)${\rm{)}}
is $1$-far from $(A_n)$.
\end{theorem}

\begin{proof}
Pick $n \in \N$, and find $f \in Y_r = \{\Re f : f \in Y\}$ such that
$\|f\| = 1 = E(f, A_n)$. Let $N = \psi(n)$. By \cite[Theorem 2.3]{voigt},
there exists $f \in Y$ and $t_1 < t_2 < \ldots < t_{N+1}$ in $[a,b]$
such that $\|f\| = 1$, and $f(t_k) = (-1)^k$ for every $k$ (the
original proof is formulated for the spaces $C(K)$, where $K$ is a
compact interval, but it can be easily extended to the general $C_0(I)$).
We have to show that $\|f-g\| \geq 1$ for any $g \in A_n \backslash \{0\}$.
As $g$ has fewer than $N$ zeroes, there exists $k \in \{1,\ldots,N\}$
such that $\Re g$ does not vanish on $(t_k, t_{k+1})$. Then $\Re g(t_k)$ and
$\Re g(t_{k+1})$ are either both non-negative, or both non-positive. Therefore,
$$
\|f-g\| \geq \max \{ |f(t_k) - \Re g(t_k)|, |f(t_{k+1}) - \Re g(t_{k+1})| \} \geq 1 ,
$$
which is what we need.
\end{proof}

\begin{corollary}\label{cor_gen_haar}
Suppose $\{A_n\}$ is a generalized Haar approximation scheme in $C[a,b]$,
with $-\infty < a < b < \infty$.
Then for all $\{\varepsilon_n\}\searrow 0$ there exists
$f\in C[a,b]$, analytic on $(a,b)$, such that
$E(f,A_n)\neq \mathbf{O}(\varepsilon_n)$.
\end{corollary}

\begin{proof}
We re-use some ideas from the proof of Corollary~\ref{suavidad}. It suffices to
consider $0 < a < b$. Then $Y=\overline{\mathbf{span}[\{x^{n^2}:n\in\N\}]}^{C[a,b]}$
is a proper subspace of $C[a,b]$, whose elements are analytic on $(a,b)$. Finally,
Theorem~\ref{gen_haar} guarantees the existence of $f \in Y$ with the desired properties.
\end{proof}

Suppose $K$ is a compact Hausdorff set. A closed subalgebra
$Y$ of $C(K)$ is called {\it uniform} if it contains the constants, and
separates points in $K$ (the disk algebra is an accessible example).

\begin{theorem}\label{uniform}
Suppose $(A_n)$ is a generalized Haar approximation scheme in $C[a,b]$,
and $Y$ is an infinite-dimensional uniform subalgebra of $C[a,b]$.
Then $Y$ {\rm{(}}equipped with the norm of $C[a,b]${\rm{)}}
is $1$-far from $(A_n)$.
\end{theorem}

\begin{proof}
Fix $N \in \N$, and $\varepsilon > 0$. 
We find distinct points $t_0, t_1, \ldots, t_N \in [a,b]$ and $h \in Y$,
such that $|h(t_j) - (-1)^j| < \varepsilon$ for every $j$, and $\|h\| < 1 + \varepsilon$.
Once this is done, pick a non-zero $f \in A_n$.
If $N > \psi(n)$, there exists $j \in \{1, \ldots, N\}$ such that
$\Re f$ doesn't change sign on $(t_{j-1}, t_j)$. Therefore,
$\max\{ |f(t_j) - (-1)^j| , |f(t_{j-1}) - (-1)^{j-1}| \geq 1$.
By the triangle inequality,
$\max\{ |f(t_j) - h(t_j)| , |f(t_{j-1}) - h(t_{j-1})|\} > 1 - \varepsilon$.
Therefore, $E(h/\|h\|, A_n) > (1-\varepsilon)/(1+\varepsilon)$.

To construct $h$, recall that $t \in [a,b]$ is a {\it peak point} if
there exists $f \in Y$ such that $|f(t)| = 1$, and $|f(s)| < 1$ for any
$s \neq t$ (we say that $f$ {\it peaks at $t$}).
The reader is referred to \cite[Section II.11]{gamelin} for more information.
In particular, the paragraph at the end of that section shows that,
for any infinite dimensional uniform algebra, the set of peak points is
infinite.

Suppose $t_0, t_1, \ldots, t_N \in [a,b]$ are peak points for
the functions $f_0, f_1, \ldots, f_N \in Y$. We can assume that
$f_j(t_j) = 1$ for every $j$. Find disjoint open sets $U_j \supset t_j$,
and $M$ so large that $|f_j(s)|^M < \varepsilon/N$ for any $s \notin U_j$
($0 \leq j \leq N$). Then $h = \sum_{j=0}^N (-1)^j f_j^M$ satisfies the
desired properties. Indeed, for $0 \leq j \leq N$,
$$
|h(t_j) - (-1)^j| \leq \sum_{k \neq j} |f_k(t_j)|^M < \varepsilon .
$$
Furthermore, for $s \in U_j$,
$$
|h(s)| \leq |f_j(s)| + \sum_{k \neq j} |f_k(s)|^M < 1 + \varepsilon ,
$$
while for $s \notin \cup_j U_j$,
$|h(s)| \leq \sum_j |f_j(s)|^M < \varepsilon$.
\end{proof}

Slightly more can be said when $A$ is the disk algebra.

\begin{theorem}\label{disk_alg}
Suppose $\varepsilon_1 > \varepsilon_2 > \ldots > 0$, and $\sum_i \varepsilon_i < \infty$.
Furthermore, suppose ${\mathcal{D}}$ is a generalized Haar system on
the unit circle $\T$. Then, for any increasing sequence $(n_i)$ of
natural numbers, there exists $f$ in the disk algebra $A$ such that
$\|f\| \leq 3 \sum_i \varepsilon_i$,
and $E(f, \Sigma_{n_i}({\mathcal{D}})) > \varepsilon_i$ for all $i$.
\end{theorem}

\begin{proof}
In the proof, we rely on Rudin-Carleson Theorem \cite[III.E.2]{Woj}:
suppose $E$ is a subset of $\T$ of measure $0$, $g$ is a continuous
function on $E$, and $h$ is a strictly positive continuous function on $\T$,
such that $h \geq |g|$ on $E$. Then there exists $f \in A$ such that
$f|_E = g$, and $|f| \leq h$.

For notational simplicity, we denote by $[t,s]$ ($t,s \in \T$) the arc,
stretching from $t$ to $s$, in the counterclockwise direction.
For each $i$, fix an even $N_i > \psi(n_i)$ (here, $\psi$ is the function
appearing in the definition of the Haar system).
We construct functions $(f_i)_{i=1}^\infty$, in such a way that,
for every $i$:
\begin{enumerate}
\item
There exists an arc $J_i$, containing points $(t_{ij})_{j=1}^{N_i}$
(enumerated counterclockwise), such that $t_{k \ell} \notin J_i$ for any $k < i$, and
$1 \leq \ell \leq N_k$.
\item
$|\sum_{k<i} \Re f_k| < \varepsilon_i/2$ on $J_i$.
\item
$\|f_i\| \leq 2 \varepsilon_i$, and
$f(t_{ij}) = 2 \varepsilon_i (-1)^j$ for $1 \leq j \leq N_i$.
\item
For $k < i$, $|f_i(t_{k\ell})| < \varepsilon_{2i}/2^{2i}$.
\end{enumerate}
Then $f = \sum_{i=1}^\infty f_i$ the desired properties. Clearly,
$\|f\| \leq 2 \sum_i \varepsilon_i$. To establish
$E(f, \Sigma_{n_i}({\mathcal{D}}))$ $> \varepsilon_i$,
pick $g \in \Sigma_{n_i}({\mathcal{D}})$. $\Re g$ cannot change sign more than
$N_i - 2$ times, hence there exists $j \in \{1, \ldots, N_i-1\}$ so that the
signs of $\Re g$ at $t_{ij}$ and $t_{i,j+1}$ are the same. Suppose $j$ is even
(the case of $j$ being odd is handled similarly). Then
$$
\Re f(t_{ij}) \geq
\Re f_i(t_{ij}) - |\sum_{k<i} \Re f_k(t_{ij})| - \sum_{k>i} |f_k(t_{ij})|
\geq 2 \varepsilon_i - \frac{\varepsilon_i}{2} - \sum_{k>i} \frac{\varepsilon_{2k}}{2^{2k}}
> \varepsilon_i ,
$$
and similarly, $\Re f(t_{i,j+1}) < - \varepsilon_i$. Therefore,
$$
\max\{|\Re f(t_{ij}) - \Re g(t_{ij})|, |\Re f(t_{i,j+1}) - \Re g(t_{i,j+1})| \} > \varepsilon_i ,
$$
which implies $\|f-g\| > \varepsilon_i$.

To define $f_1$, fix an even $N_1 > \psi(n_1)$, and select points
$(t_{1i})_{i=1}^{N_1}$ (enumerated counterclockwise). By Rudin-Carleson Theorem,
there exists $f_1 \in A$ such that $f_1(t_i) = 2 (-1)^i \varepsilon_1$ for any
$1 \leq i \leq N_1$.

Now suppose $f_1, \ldots, f_{i-1}$, and the points $t_{k \ell}$ ($1 \leq k < i$,
$1 \leq \ell \leq N_k$) with the desired properties have already been defined.
Then, for $1 \leq \ell \leq N_1$,
$$
\big| 2 \varepsilon_1 (-1)^\ell - \sum_{k=1}^{i-1} f_k (t_{1j}) \big| \leq
\sum_{k=2}^{i-1} |f_k (t_{1j})| \leq
\sum_{k=2}^{i-1} \frac{\varepsilon_{2k}}{2^{2k}} < \frac{3 \varepsilon_1}{2} .
$$
In particular, $\Re(\sum_{k=1}^{i-1} f_k)$ changes sign at least $N_1$ times.
Use the continuity of $f_1, \ldots, f_{k-1}$ to find an arc $J_i$, not containing
any of the points $t_{k \ell}$ ($k < i$, $1 \leq \ell \leq N_k$),
and such that $|\sum_{k<i} \Re f_k| < \varepsilon_i/2$ on $J_i$.
Rudin-Carleson Theorem guarantees the existence of a function $f_i$
with desired properties.
\end{proof}

Recall that a sequence $(e_i)_{i=1}^{\infty}\subset X$ in a Banach space $X$
is called a {\it basis} in $X$ if for e\-ve\-ry $x\in X$ there exists a unique sequence of
scalars $(a_n(x))$ such that $x=\sum_{n=1}^{\infty}a_n(x)e_n$.
In this case, the {\it basis projections} $P_n$, defined by $P_n(x)=\sum_{k=1}^na_{k}(x)e_k$,
are uniformly bounded. The basis $(e_i)$ is called {\it $C$-unconditional} if,
for any finite sequences of scalars $(a_i)$ and $(b_i)$,
$\|\sum a_i b_i e_i\| \leq C (\sup_i |b_i|) \|\sum a_k e_k\|$.
A basis is {\it unconditional} if it is $C$-unconditional for some $C$.
It is easy to see that every Banach space with an unconditional basis can be
renormed to make this basis $1$-unconditional.
We refer the reader to \cite{AK} or \cite{LZ} for more information about bases.


\begin{theorem}\label{unc_basis}
Suppose $Y$ is a closed infinite dimensional subspace of a Banach space $X$.
and $(e_i)_{i \in \N}$ is an unconditional basis in $X$.
Then $Y$ satisfies Shapiro's Theorem
with respect to $(X, \{\Sigma_n(\{e_1, e_2, \ldots\})\})$.
\end{theorem}

The condition of $Y$ being closed (in the norm inherited from $X$)
cannot be omitted, even when $\overline{Y} = X$. Consider, for instance,
the canonical basis $(e_i)$ in $X = \ell_2$, and let $Y$ be the set of all
$x = (x_1, x_2, \ldots) \in X$, for which $\|x\|_Y = (\sum_k k^2 |x_k|^2)^{1/2}$
is finite. Then $\overline{Y} = X$. However, for all $x\in Y$,
$$
E(x, \Sigma_n(\{e_1, e_2, \ldots\}))^2 \leq \sum_{k=n+1}^\infty |x_k|^2
\leq (n+1)^{-2} \sum_{k=n+1}^\infty k^2 |x_k|^2 \leq (n+1)^{-2} \|x\|_Y ,
$$
hence $E(x, \Sigma_n(\{e_1, e_2, \ldots\})) = \mathbf{O}(n^{-1})$.

To prove Theorem~\ref{unc_basis}, it suffices to combine Theorem~\ref{far_shap} with

\begin{theorem}\label{uncond_far}
Suppose $(e_i)$ is a $1$-unconditional basis in a Banach space $X$.
Then any closed infinite dimensional subspace of $X$ is $1$-far from
the approximation scheme $\{\Sigma_n(\{e_1, e_2, \ldots\})\}$.
\end{theorem}

\begin{proof}
Fix $n \in \N$ and $c > 1$. We have to show that there exists
$x \in Y$ such that
\begin{equation}
\|x\| < c \, \, \, {\mathrm{and}} \, \, \,
E(x, \Sigma_n(\{e_1, e_2, \ldots\})) > 1/c .
\label{requirement}
\end{equation}
To this end, pick
$\sigma \in (0,1)$, $\delta \in (0,\sigma/2)$, and $M > n$ in such a way that
$$
(1+\sigma)^2 < c , \, \, \,
\frac{1}{1+\sigma} - \sigma (1+\delta) > \frac{1}{c} , \, \, \,
{\mathrm{and}} \, \, \,
\frac{M-n}{M} > \frac{1+\delta}{1+\sigma} .
$$

Recall that the basis projections $P_n$ are defined by setting
$P_n(\sum_i a_i e_i) = \sum_{i \leq n} a_i e_i$. For notational convenience,
we put $P_0 = 0$. As the basis $(e_i)$ is $1$-unconditional, the projections
$P_n$ and $I - P_n$ are contractive for every $n$.

Find $0 = N_0 < N_1 < N_2 < \ldots$ and
$y_1, y_2 \ldots \in Y$ so that, for every $i$, $\|y_i\| = 1$,
$(I - Q_{i-1}) y_i = 0$, and $\|Q_j y_i\| < \delta 4^{-j}$
for $j \geq i$ (here, $Q_j = I - P_{N_j}$).
Indeed, suppose $N_0 < \ldots < N_k$ and $y_1, \ldots, y_k$ ($k \geq 0$)
with desired properties have already been selected. Then
$X^{(k)} = \{x \in X : Q_k x = 0\}$ is a finite codimensional subspace of $X$,
hence $Y \cap X^{(k)}$ is non-empty. Pick $y_{k+1} \in Y \cap X^{(k)}$ of norm $1$.
Note that $\lim_m \|(I - P_m) x\| = 0$ for any $x \in X$, hence we can find
$N_{k+1} > N_k$ such that $\|Q_{k+1} y_i\| < \delta 4^{-(k+1)}$ for $1 \leq i \leq k+1$.

Let $y_i^\prime = y_i - Q_i y_i$. Then $\|y_i - y_i^\prime\| < \delta/4^i$
for every $i \in \N$, hence $1 - \delta/4^i < \|y_i^\prime\| \leq 1$.
Furthermore, the vectors $y_i^\prime$ have disjoint support:
$y_i^\prime \in \spn[e_s : N_{i-1} < s \leq N_i]$. Therefore,
\begin{equation}
\sum_i |\alpha_i| \geq \|\sum_i \alpha_i y_i^\prime\| \geq
\max_i (1 - \delta/4^i) |\alpha_i| \geq \frac{1}{2} \max_i |\alpha_i|
\label{disj_supp}
\end{equation}
for any finite sequence $(\alpha_i)$. Consider a linear map
$T : \spn[y_i^\prime : i \in \N] \to \spn[y_i : i \in \N]$, defined by
$T y_i^\prime = y_i$. Then
\begin{equation}
\|T - I\| < \sigma . 
\label{norm_T}
\end{equation}
Indeed, suppose $\|\sum_i \alpha_i y_i^\prime\| = 1$. By \eqref{disj_supp},
$$
\|(T - I)(\sum_i \alpha_i y_i^\prime)\| \leq
\sum_i |\alpha_i| \|y_i - y_i^\prime\| \leq
2 \delta \sum_i 4^{-i} < \sigma .
$$

By Krivine's Theorem (see e.g.~\cite[Section 11.3]{AK}), there exists
$q \in [1,\infty]$, $1 < p_0 < \ldots < p_M$, and norm $1$ vectors
$z_j^\prime = \sum_{i=p_{i-1}}^{p_i-1} \beta_i y_i^\prime$, such that
\begin{equation}
\frac{1}{1+\delta} \Big( \sum_{j=1}^M |\gamma_j|^q \Big)^{1/q} \leq
\|\sum_{j=1}^M \gamma_j z_j^\prime\| \leq
(1+\delta) \Big( \sum_{j=1}^M |\gamma_j|^q \Big)^{1/q} .
\label{krivine}
\end{equation}
Let $x^\prime = M^{-1/q} \sum_{j=1}^M z_j^\prime$. By the above,
$\|x^\prime\| \leq 1 + \delta$.
We claim that, for any sequence $(\alpha_i)$ with at most $n$ non-zero entries,
$\|x^\prime - \sum_i \alpha_i e_i\| > 1/(1+\sigma)$.
Indeed, let $S$ be the set of $j \in \{1, \ldots, M\}$ with the
property that $\alpha_i = 0$ whenever $p_{j-1} \leq i < p_j$.
Define the projection $R$ by setting $R e_i = e_i$ if
$i \in \cup_{j \in S} [p_{j-1}, p_j)$, $R e_i = 0$ otherwise.
By the $1$-unconditionality of $(e_i)$, the projection $R$ is contractive.
By \eqref{krivine},
$$
\|x^\prime - \sum_i \alpha_i e_i\| \geq
\|R(x^\prime - \sum_i \alpha_i e_i)\| =
M^{-1/q} \|\sum_{j \in S} z_j^\prime\| \geq
\frac{1}{1+\delta} \Big( \frac{M-n}{M} \Big)^{1/q} > \frac{1}{1+\sigma} .
$$


It remains to show that $x = T x^\prime$ satisfies \eqref{requirement}.
By \eqref{norm_T}, $\|x - x^\prime\| \leq \|T - I\| \|x^\prime\| < \sigma (1+\delta)$,
and therefore, by the triangle inequality,
$\|x\| \leq \|x^\prime (1+\sigma)(1+\delta) < c$.
Furthermore, if a sequence $(\alpha_i)$ with at most $n$ non-zero entries, then
$$
\|x - \sum_i \alpha_i e_i\| \geq
\|x^\prime - \sum_i \alpha_i e_i\| - \|x - x^\prime\| \geq
\frac{1}{1+\sigma} - \sigma (1+\delta) > \frac{1}{c} ,
$$
which yields \eqref{requirement}. 
\end{proof}

Finally, we deal with non-commutative sequence spaces.
Suppose $\se$ is a separable symmetric sequence space.
Consider the approximation scheme $(A_i)$ in $\ce$, where $A_i$ is
the space of operators of rank not exceeding $i$. Reasoning as in
\cite[Section 6.5]{almira_oikhberg}, we see that $(\ce, \{A_i\})$ satisfies
Shapiro's Theorem. A stronger statement holds.

\begin{proposition}\label{operators}
Suppose $\se$ is a symmetric sequence space, and the approximation scheme
$(A_i)$ is defined as above. Then every finite codimensional subspace of
$\ce$ is $1$-far from $(A_i)$.
\end{proposition}

\begin{proof}
Let $(e_i)$ be an orthonormal basis in $\ell_2$. Denote by $Z$ the
space of operators $T \in \ce$ which are diagonal relative to $(e_i)$
(that is, $T e_i = s_i e_i$ for every $i$). Define a map $U : \se \to Z$, taking
$s = (s_1, s_2, \ldots)$ to the operator $U(s)$, defined via $U(s) e_i = s_i e_i$.
Clearly, $U$ is an isometry.
Denote the canonical basis of $\se$ by $(\delta_i)$, and let
$B_i = \Sigma_i(\{\delta_1, \delta_2, \ldots\})$. By
\eqref{best_appr}, $E(s, B_i) = E(U(s), A_i)$ for every $s$ and $i$.

Note that $Y^\prime = Y \cap Z$ is a finite codimensional subspace of $Z$.
By Proposition~\ref{uncond_far}, $E(S(U^{-1}(Y^\prime)), B_i) = 1$ for
every $i$. Therefore, $E(S(Y^\prime), A_i) = 1$ for every $i$.
\end{proof}

\section{Acknowledgement}
We are very grateful to an anonymous referee of this paper for his/her many interesting comments. His/her advice has greatly helped to improve its readability.


\bigskip

\footnotesize{J. M. Almira

Departamento de Matem\'{a}ticas. Universidad de Ja\'{e}n.

E.P.S. Linares,  C/Alfonso X el Sabio, 28

23700 Linares (Ja\'{e}n) Spain

email: jmalmira@ujaen.es}

\bigskip

\footnotesize{T. Oikhberg

Department of Mathematics, The University of California at Irvine, Irvine CA 92697, {\it and}

Department of Mathematics, University of Illinois at Urbana-Champaign, Urbana, IL 61801

email: toikhber@math.uci.edu}

\end{document}